\documentclass[reqno, 11pt]{amsart}

\usepackage[text={155mm,235mm},centering]{geometry}
\input diagxy
\input xy

\usepackage[active]{srcltx} 

\newtheorem{thm}{Theorem}[section]
\newtheorem{cor}[thm]{Corollary}
\newtheorem{lem}[thm]{Lemma}
\newtheorem{prop}[thm]{Proposition}

\theoremstyle{definition}
\newtheorem{defn}[thm]{Definition}
\theoremstyle{property}

\theoremstyle{remark}
\newtheorem{rem}[thm]{Remark}

\newtheorem{main theorem}[thm]{Main Theorem}

\newtheorem{ex}[thm]{Example}
\numberwithin{equation}{section}
\usepackage[mathscr]{eucal}
\usepackage[all]{xy}
\usepackage{mathrsfs}
\usepackage{xypic}
\usepackage{amsfonts}
\usepackage{amsmath}
\usepackage{amsthm,bm}
\usepackage{amssymb,tikz-cd}
\usepackage{latexsym}
\usepackage{tabularx}
\usepackage{graphicx}
\usepackage{pict2e}
\usepackage[pagebackref]{hyperref}
\usepackage{tikz}
\usepackage{textcomp}
\usepackage[scr=rsfs]{mathalfa}
\definecolor{ceruleanblue}{rgb}{0.16, 0.32, 0.75}
\hypersetup{colorlinks=true,allcolors=ceruleanblue}
\allowdisplaybreaks

\begin{document}

\title[Morse-Novikov cohomology for blow-ups]{Morse-Novikov cohomology for  blow-ups of complex manifolds}

\author{Lingxu Meng}
\address{Department of Mathematics, North University of China, Taiyuan, Shanxi 030051, P.R. China}
\email{menglingxu@nuc.edu.cn}%

\subjclass[2010]{Primary 53C56; Secondary 55N35, 32Q55}
\keywords{Morse-Novikov cohomology; blow-up; weight $\theta$-sheaf; self-intersection formula; Leray-Hirsch theorem}


\begin{abstract}
  The weight $\theta$-sheaf $\underline{\mathbb{R}}_{X,\theta}$ helps us to reinterpret  Morse-Novikov cohomologies via sheaf theory.
  We give several theorems of K\"{u}nneth  and  Leray-Hirsch types.
  As applications, we prove that the $\theta$-Lefschetz number is independent of $\theta$ and calculate the Morse-Novikov cohomologies of projective bundles.
  Based on these results, we give two blow-up formulae  on (\emph{not necessarily  compact}) complex manifolds, where the self-intersection formulae play a key role in  establishing the explicit expressions for them.

\end{abstract}

\maketitle

\section{Introduction}
In this paper, all smooth manifolds are assumed to be \emph{connected}, \emph{paracompact} and all submanifolds (resp. complex submanifolds) are  assumed to be \emph{closed} (in the topological sense) \emph{embedded smooth} (resp. \emph{complex}) \emph{submanifolds   without boundary}.
For a smooth manifold $X$ and a closed one-form $\theta$ on $X$,  let $\mathcal{A}^p(X)$ be the space of smooth $p$-forms and  define $d_\theta:\mathcal{A}^p(X)\rightarrow\mathcal{A}^{p+1}(X)$ as $d_\theta\alpha=d\alpha+\theta\wedge\alpha$ for any $\alpha\in\mathcal{A}^p(X)$.
Obviously, $d_\theta\circ d_\theta=0$, so $(\mathcal{A}^\bullet(X), d_\theta)$ is a complex.
Denote its $p$-th cohomology by $H^p_\theta(X)$, which is called the \emph{Morse-Novikov cohomology} \cite{Mi,OV,YZ}, \emph{Lichnerowicz cohomology} \cite{BK,LLMP},  \emph{adapted cohomology} \cite{DO,Vai}, or \emph{$d_\theta$-cohomology} \cite{Ba1,Ba2,HR}.
In this article, we call it the \emph{Morse-Novikov cohomology}.
Similarly,  $H^p_{\theta,c}(X)$ is defined as the $p$-th cohomology group of the complex $(\mathcal{A}_c^\bullet(X), d_\theta)$ of the spaces of smooth forms with compact supports, which are called the $p$-th  \emph{Morse-Novikov cohomology with compact support}. Clearly, if $\theta=0$, they are just the de Rham cohomology and the de Rham cohomology with compact support.

The Morse-Novikov cohomology was originally defined by A. Lichnerowicz  \cite{L} and D. Sullivan \cite{S} in the context of Poisson geometry and infinitesimal computations in topology, respectively.
It was well used to study the locally conformally K\"{a}hlerian and locally conformally symplectic  structures \cite{Ba1,Ba2,BK,DO,HR,LLMP,Vai}.
S. Novikov \cite{N} introduced a generalization of the classical Morse theory to the case of circle-valued Morse functions.
A. Pajitnov \cite{P1}  observed the relation of the circle-valued Morse theory to the homology with local coefficients and the perturbed de Rham differential; see also \cite[p. 414-416]{P2}.

Unfortunately, the Morse-Novikov cohomology is much more difficult to calculate than the de Rham cohomology, since it depends on a closed one-form.
Until now, we don't know much  about the Morse-Novikov cohomology.
For instance,
the Morse-Novikov cohomology of compact manifolds is  finitely dimensional  \cite{Vai} and
the Mayer-Vietoris sequence, Poincar\'{e} duality theorem hold for the Morse-Novikov cohomology \cite{HR}.
For completely solvable solvmanifolds, the Morse-Novikov cohomology coincides with the cohomology of the corresponding Lie algebra   \cite{Mi}.
A well-known result is that a compact Riemannian manifold $X$ endowed with a parallel one-form $\theta$ has trivial  Morse-Novikov cohomology \cite{LLMP,OV}.
By Atiyah-Singer index theorem,  the Euler characteristic of Morse-Novikov cohomology coincides with the classical Euler characteristic  \cite{BK}.
For a Morse's function $f$ and $\lambda\in\mathbb{R}$, $\textrm{d}_{\lambda \textrm{d}f}$ is just the Witten deformation, which was used to prove  strong Morse inequalities by E. Witten \cite{W}.
Moreover, there are some results for specific examples, see \cite{Ba2,HR,Mi,O1,O2}.
More geometric, topological  and dynamical applications of Morse-Novikov cohomology can be found in  \cite{F,P2}.

On complex manifolds, if $\theta$ is a complex closed one-form,  the cohomology $H_{\theta}^*(X)$ can be viewed as the cohomology of  a flat bundle, i.e., the weight line bundle \cite{N,OV,YZ}, or a locally constant sheaf of $\mathbb{C}$-modules of rank $1$ \cite{S}.
As we know, the two viewpoints are equivalent.
We will study the Morse-Novikov cohomologies by the language of locally constant sheaves, which is much more convenient.

The self-intersection formula is important in the intersection theory.
We establish it for Morse-Novikov cohomologies as follows, which plays a key role in writing out the  explicit expressions of blow-up formulae.
\begin{lem}[Self-intersection formulae]\label{key}
Let $Y$ be an oriented  submanifold of an oriented  smooth manifold $X$ of codimension $r$.
Denote by $i:Y\rightarrow X$  the inclusion and by $[Y]\in H_{dR}^r(X)$  the fundamental class of $Y$ in $X$.
Assume that $\theta$ is a closed one-form on $X$.
Then
\begin{displaymath}
i^*i_*\sigma=[Y]|_Y\cup\sigma.
\end{displaymath}
for $\sigma\in H_{\theta|_Y}^*(Y)$ or $H_{\theta|_Y,c}^*(Y)$.
\end{lem}

X.-D Yang and G. Zhao \cite{YZ} proved that there exists a  blow-up formula of Morse-Novikov cohomology under some assumptions on closed one-forms for  compact locally conformally K\"{a}hlerian manifolds, which generalized the corresponding result of singular cohomology for compact K\"{a}hler manifolds \cite[Proposition 13.1]{Le}\cite[Theorem 7.31]{Vo}.
It  seems difficult to write out the expression explicitly.
We will establish blow-up formulae of Morse-Novikov cohomologies without additional assumptions for arbitrary complex manifolds.
Moreover, we express them explicitly by Theorem \ref{L-H} and Lemma \ref{key} as follows.
\begin{thm}[Blow-up formulae]\label{1.3}
Let $\pi:\widetilde{X}\rightarrow X$ be the blow-up of a complex manifold $X$ along a complex submanifold $Y$ of complex codimension $r$.
Assume that $\theta$ is a  closed one-form on $X$ and $\tilde{\theta}=\pi^*\theta$.
Set $E=\pi^{-1}(Y)$ and let $i_E:E\rightarrow \widetilde{X}$ be the inclusion.
Then
\begin{equation}\label{b-u-m}
\pi^*+\sum_{i=1}^{r-1}i_{E*}\circ (h^{i-1}\cup)\circ (\pi|_E)^*
\end{equation}
gives  isomorphisms
\begin{equation}\label{isomorphism 1}
H_{\theta}^k(X)\oplus \bigoplus_{i=1}^{r-1}H_{\theta|_Y}^{k-2i}(Y)\tilde{\rightarrow} H_{\tilde{\theta}}^k(\widetilde{X}),
\end{equation}
\begin{equation}\label{isomorphism 2}
H_{\theta,c}^k(X)\oplus \bigoplus_{i=1}^{r-1}H_{\theta|_Y,c}^{k-2i}(Y)\tilde{\rightarrow} H_{\tilde{\theta},c}^k(\widetilde{X})
\end{equation}
for any $k$, where $\pi|_E:E\rightarrow Y$ is viewed as the complex projectivization $E=\mathbb{P}(N_{Y/X})$ of the normal bundle $N_{Y/X}$ of $Y$ in $X$ and $h=c_1(\mathcal{O}_E(-1))\in H_{dR}^2(E)$ is the first Chern class of the universal line bundle $\mathcal{O}_E(-1)$ over $E$.
\end{thm}


\begin{rem}
For a smooth  manifold $X$, set $\mathcal{A}_{X,\mathbb{C}}^k=\mathcal{A}_X^k\otimes_{\mathbb{R}}\mathbb{C}$ and $\mathcal{D}_{X,\mathbb{C}}^{'k}=\mathcal{D}_X^{'k}\otimes_{\mathbb{R}}\mathbb{C}$.
For a complex closed one-form $\theta$, we can define the  Morse-Novikov cohomology 
as the one in the real case. \emph{All arguments in this article hold for both real and complex cases}.
In what follows, we only consider the cases for  \emph{real} closed one-forms.
\end{rem}

Using the expression (\ref{b-u-m}), we can easily generalize Theorem \ref{1.3} to the cases of the cohomologies with values in local systems by Mayer-Vietoris sequences \cite{Me2}.
After the first version  \cite{Me1} of the present paper,  Y. Chen, S. Yang \cite{CY} and Y. Zou \cite{Z} used different proofs to obtain the blow-up formula  on compact complex manifolds without  proofs of self-intersection formula and expression (\ref{b-u-m}). They depend more on sheaf theory (more precisely says, the local systems) and  the  compactness is necessary there.
The present paper contains  many interesting results and  inimitable techniques for  Morse-Novikov cohomology, which seem difficult to be generalized  to the cohomology  with the value in a general local system (see Sect. 6 for more explanations) and are valuable for the study of locally conformally K\"{a}hlerian and locally conformally symplectic structures.

The paper is organized as follows.
In Sect. 2, we study the weight $\theta$-sheaf $\underline{\mathbb{R}}_{X,\theta}$ and its relationship with the Morse-Novikov cohomology.
In Sect. 3, the K\"{u}nneth theorems  are established for Morse-Novikov cohomologies.
As an application, we  generalize a result of  G. Bande and D. Kotschick on the Euler characteristic with a much more elementary proof.
In Sects. 4-5, Theorem \ref{L-H}, Lemma \ref{key} and Theorem \ref{1.3} are proved respectively.
In Sect. 6, we explain the particularity of the Morse-Novikov cohomology compared with the cohomology  with the value in a general local system. \\\\
\textbf{Notations.} We fix some notations in this article.

- $H^*$\mbox{ } the graded vector space $\bigoplus\limits_{p\geq0}H^p$;

- $H_1^*\otimes H_2^*$\mbox{ }  the graded vector space satisfies that  $(H_1^*\otimes H_2^*)^p=\bigoplus\limits_{r+s=p}H_1^r\otimes H_2^s$.\\
Assume that $X$ is a smooth manifold and $Y$ is a closed submanifold of $X$.

- $\mathcal{A}^p(X)$ (resp. $\mathcal{A}_c^p(X)$) the space of real-valued smooth $p$-forms (resp. smooth $p$-forms with compact supports) on $X$;

- $\mathcal{A}_X^p$ the sheaf of germs of real-valued smooth $p$-forms on $X$;

- $\textrm{dim}X$ the (real) dimension of $X$;

- $\textrm{codim}Y$ the (real) codimension of $Y$ in $X$.\\
In addition, assume that $X$ is oriented.

- $\mathcal{D}^{\prime p}(X)$ (resp. $\mathcal{D}_c^{\prime p}(X)$) the space of real-valued  $p$-currents (resp.  $p$-currents with compact supports)  on $X$;

- $\mathcal{D}_X^{\prime p}$  the sheaf of germs of $p$-currents on $X$.\\

\section{Morse-Novikov cohomology}
\subsection{Weight $\theta$-sheaf}
Let $X$ be an  $n$-dimensional smooth manifold and  $\theta$ a closed one-form on $X$.
For an open subset $U\subseteq X$, define $d_{\theta|_U}:\mathcal{A}^p(U)\rightarrow\mathcal{A}^{p+1}(U)$ as $d_{\theta|_U}\alpha=d\alpha+\theta|_U\wedge\alpha$ for $\alpha\in\mathcal{A}^p(U)$.
One easily checks that all $d_{\theta|_U}$ for open subsets $U\subseteq X$  give a morphism $d_\theta:\mathcal{A}_X^p\rightarrow\mathcal{A}_X^{p+1}$ of sheaves of $\mathbb{R}$-modules.
Clearly, $d_\theta\circ d_\theta=0$.

\begin{defn}
The kernel of $d_\theta:\mathcal{A}_X^0\rightarrow\mathcal{A}_X^1$ is called the \emph{weight $\theta$-sheaf}, denoted by $\underline{\mathbb{R}}_{X,\theta}$.
\end{defn}
Locally, $\theta=du$ for a smooth function $u$, so $d_\theta=e^{-u}\circ d\circ e^u$ and $\underline{\mathbb{R}}_{X,\theta}=\mathbb{R}e^{-u}$. Hence, the weight $\theta$-sheaf $\underline{\mathbb{R}}_{X,\theta}$ is a locally constant sheaf of $\mathbb{R}$-modules of rank $1$  and there is a soft resolution of  $\underline{\mathbb{R}}_{X,\theta}$
\begin{displaymath}
\xymatrix{
0\ar[r] &\underline{\mathbb{R}}_{X,\theta}\ar[r]^{i} &\mathcal{A}_X^0\ar[r]^{d_\theta} &\mathcal{A}_X^1\ar[r]^{d_\theta}&\cdots\ar[r]^{d_\theta}&\mathcal{A}_X^n\ar[r]&0,
}
\end{displaymath}
where $i$ is the inclusion.

Assume that $X$ is \emph{oriented}.
For any open subset $U\subseteq X$, define $d_{\theta|_U}:\mathcal{D}^{\prime p}(U)\rightarrow\mathcal{D}^{\prime p+1}(U)$ as $d_{\theta|_U}T=dT+\theta|_U\wedge T$ for $T\in\mathcal{D}^{\prime p}(U)$.
One easily checks that all $d_{\theta|_U}$ for open subsets $U\subseteq X$  give a morphism $d_\theta:\mathcal{D}_X^{\prime p}\rightarrow\mathcal{D}_X^{\prime p+1}$ of sheaves of $\mathbb{R}$-modules.
For any $T\in \mathcal{D}^{\prime p}(U)$ and $\alpha\in\mathcal{A}_c^{n-p-1}(U)$, $d_{\theta|_U}T(\alpha)=(-1)^{p+1}T(d_{-\theta|_U}\alpha)$,
so $d_\theta\circ d_\theta=0$.
Locally, $d_\theta=e^{-u}\circ d\circ e^u$ on $\mathcal{D}_X^{\prime*}$ for a smooth function $u$, which implies  $\underline{\mathbb{R}}_{X,\theta}=\textrm{ker}(d_\theta:\mathcal{D}_X^{\prime 0}\rightarrow\mathcal{D}_X^{\prime 1})$.
There is another soft  resolution of  $\underline{\mathbb{R}}_{X,\theta}$
\begin{displaymath}
\xymatrix{
0\ar[r] &\underline{\mathbb{R}}_{X,\theta}\ar[r]^{i} &\mathcal{D}_X^{\prime 0}\ar[r]^{d_\theta} &\mathcal{D}_X^{\prime 1}\ar[r]^{d_\theta}&\cdots\ar[r]^{d_\theta}&\mathcal{D}_X^{\prime n}\ar[r]&0,
}
\end{displaymath}
where $i$ is the inclusion.

\begin{lem}\label{lem fun}
Let $X$ be a  smooth manifold and $\theta$ a closed one-form on $X$. Denote by $\underline{\mathbb{R}}_{X}$  the constant sheaf with stalk $\mathbb{R}$ on $X$.

$(1)$  There exists an isomorphism $\underline{\mathbb{R}}_{X,\theta}\cong\underline{\mathbb{R}}_{X}$ if and only if $\theta$ is exact. More precisely, if $\theta=du$ for $u\in \mathcal{A}^0(X)$, then $h\mapsto e^{u}\cdot h$ gives an isomorphism $\underline{\mathbb{R}}_{X,\theta}\tilde{\rightarrow}\underline{\mathbb{R}}_X$ of sheaves.

$(2)$ For a closed one-form $\mu$ on $X$,  the tensor product $\underline{\mathbb{R}}_{X,\theta}\otimes_{\underline{\mathbb{R}}_X}\underline{\mathbb{R}}_{X,\mu}\cong\underline{\mathbb{R}}_{X,\theta+\mu }$. In particular, the dual sheaf $\underline{\mathbb{R}}^\vee_{X,\theta}\cong\underline{\mathbb{R}}_{X,-\theta}$.

$(3)$ Suppose that $f:Y\rightarrow X$ is a smooth map between smooth manifolds.
Then the inverse image sheaf $f^{-1}\underline{\mathbb{R}}_{X,\theta}\cong\underline{\mathbb{R}}_{Y,f^*\theta}$.

$(4)$ Suppose that $Y$ is a smooth manifold and $\eta$ is a closed one-form on $Y$. Let $pr_1$ and  $pr_2$ be projections of $X\times Y$ onto $X$ and $Y$ respectively. Then the external product $\underline{\mathbb{R}}_{X,\theta}\boxtimes \underline{\mathbb{R}}_{Y,\eta}\cong \underline{\mathbb{R}}_{X\times Y,pr_1^*\theta+pr_2^*\eta}$.
\end{lem}

\begin{proof}
Assume that $\underline{\mathbb{R}}_{X,\theta}$ is a constant sheaf.
Then $\{f\in \mathcal{A}^0(X)|\mbox{ }d_\theta f=0\}=\Gamma(X,\underline{\mathbb{R}}_{X,\theta})\cong\mathbb{R}$.
By solving a simple first-order ordinary differential equation,  $\theta$ is exact on $X$, see \cite[Example 1.6]{HR} for details.
Conversely, if $\theta=du$, $\underline{\mathbb{R}}_{X,\theta}=\mathbb{R}e^{-u}$.
We get  $(1)$.
Locally, $\theta=du$ and $\mu=dv$ for smooth functions $u$ and $v$.
Then  $\underline{\mathbb{R}}_{X,\theta}=\mathbb{R}e^{-u}$, $\underline{\mathbb{R}}_{X,\mu}=\mathbb{R}e^{-v}$ and $\underline{\mathbb{R}}_{X,\theta+\mu}=\mathbb{R}e^{-u-v}$ locally.
Evidently,  products of functions give an isomorphism $\underline{\mathbb{R}}_{X,\theta}\otimes_{\underline{\mathbb{R}}_X}\underline{\mathbb{R}}_{X,\mu}\tilde{\rightarrow}\underline{\mathbb{R}}_{X,\theta+\mu }$ of sheaves, i.e., $(2)$ holds.
Locally, $\theta=du$ for smooth functions $u$ and then $\underline{\mathbb{R}}_{X,\theta}=\mathbb{R}e^{-u}$,  $\underline{\mathbb{R}}_{Y,f^*\theta}=\mathbb{R}e^{-f^*u}$.
So the pullbacks of functions give an isomorphism $f^{-1}\underline{\mathbb{R}}_{X,\theta}\tilde{\rightarrow} \underline{\mathbb{R}}_{Y,f^*\theta}$, i.e., $(3)$ holds.
By $(2)$ and $(3)$, we get $(4)$ immediately.
\end{proof}

\subsection{Morse-Novikov cohomology}
Suppose that $\theta$ is a closed one-form on a smooth manifold $X$.
The cohomologies $H^*_{\theta}(X)=H^*(\mathcal{A}^\bullet(X),d_\theta)$ and  $H^*_{\theta,c}(X)=H^*(\mathcal{A}_c^\bullet(X),d_\theta)$  are said to be the  \emph{Morse-Novikov cohomology} and \emph{Morse-Novikov cohomology with compact supports} respectively.
Let $(\mathcal{A}_X^\bullet,d_\theta)\rightarrow \mathcal{I}^\bullet$ be an injective resolution of the complex $(\mathcal{A}_X^\bullet,d_\theta)$ of sheaves in the category of sheaves on $X$.
Then it induces   isomorphisms
\begin{displaymath}
H_{\theta}^*(X)= H^*(\mathcal{A}^\bullet(X),d_{\theta})\tilde{\rightarrow} H^*(\Gamma(X,\mathcal{I}^\bullet))=H^*(X,\underline{\mathbb{R}}_{X,\theta}),
\end{displaymath}
\begin{displaymath}
H_{\theta,c}^*(X)= H^*(\mathcal{A}_c^\bullet(X),d_{\theta})\tilde{\rightarrow} H^*(\Gamma_c(X,\mathcal{I}^\bullet))=H_c^*(X,\underline{\mathbb{R}}_{X,\theta}),
\end{displaymath}
which are both denoted by $\rho$.
That is to say, \emph{the two kinds of Morse-Novikov cohomologies  can be both viewed as the cohomologies  of the weight $\theta$-sheaf $\underline{\mathbb{R}}_{X,\theta}$ via $\rho$.}
For a $d_\theta$-closed $\alpha\in\mathcal{A}^*(X)$  (resp. $\mathcal{A}_c^*(X)$), denote by $[\alpha]_{\theta}$ (resp. $[\alpha]_{\theta,c}$) its class in $H_\theta^*(X)$ (resp. $H_{\theta,c}^*(X)$).
Moreover, assume that $X$ is \emph{oriented}.
The natural inclusion $(\mathcal{A}_X^\bullet,d_\theta)\hookrightarrow(\mathcal{D}_X^{\prime\bullet},d_\theta)$ induces an isomorphism $H_{\theta}^*(X)\tilde{\rightarrow} H^*(\mathcal{D}^{\prime\bullet}(X),d_{\theta})$ and $H_{\theta,c}^*(X)\tilde{\rightarrow} H^*(\mathcal{D}_c^{\prime\bullet}(X),d_{\theta})$. 
We will \emph{not distinguish} $H_{\theta}^*(X)$ and $H^*(\mathcal{D}^{\prime\bullet}(X),d_{\theta})$ (resp. $H_{\theta,c}^*(X)$ and $H^*(\mathcal{D}_c^{\prime\bullet}(X),d_{\theta})$).
For a $d_\theta$-closed $T\in\mathcal{D}^{\prime*}(X)$ (resp. $\mathcal{D}_c^{\prime*}(X)$), denote by $[T]_{\theta}$ (resp. $[T]_{\theta,c}$)  its class in $H_{\theta}^*(X)$ (resp. $H_{\theta,c}^*(X)$).
Assume that $u$  is  a smooth function  on $X$.
The isomorphism $e^{-u}\cdot:(\mathcal{A}_X^\bullet,d_\theta)\rightarrow(\mathcal{A}_X^\bullet,d_{\theta+\textrm{d}u})$ of complexes of sheaves induces  isomorphisms
\begin{equation}\label{isom}
e^{-u}\cdot:H_{\theta}^*(X)\rightarrow H_{\theta+du}^*(X),\mbox{ }[\alpha]_{\theta}\mapsto[e^{-u}\alpha]_{\theta+du},
\end{equation}
\begin{equation}\label{isom2}
e^{-u}\cdot:H_{\theta,c}^*(X)\rightarrow H_{\theta+du,c}^*(X),\mbox{ }[\alpha]_{\theta,c}\mapsto[e^{-u}\alpha]_{\theta+du,c}.
\end{equation}

For an open set $W\subseteq X$, we briefly write $H^*_{\theta|_W}(W)$ and $H^*_{\theta|_W,c}(W)$ as $H^*_{\theta}(W)$ and $H^*_{\theta,c}(W)$ respectively.

\subsection{Pullback and Pushforward}
Let $f:X\rightarrow Y$ be a smooth map between oriented smooth manifolds and $\theta$ a closed one-form on $Y$. Set $\tilde{\theta}=f^*\theta$ and $r=\textrm{dim}X-\textrm{dim}Y$. Define the \emph{pullback} $f^*:H^*_{\theta}(Y)\rightarrow H^*_{\tilde{\theta}}(X)$ as $[\alpha]_{\theta}\mapsto [f^*\alpha]_{\tilde{\theta}}$ and the \emph{pushforward} $f_*:H^*_{\tilde{\theta},c}(X)\rightarrow H^{*-r}_{\theta,c}(Y)$ as $[T]_{\theta,c}\mapsto [f_*T]_{\tilde{\theta},c}$. Moreover, if $f$ is proper, we can also define $f^*:H^*_{\theta,c}(Y)\rightarrow H^*_{\tilde{\theta},c}(X)$ and $f_*:H^*_{\tilde{\theta}}(X)\rightarrow H^{*-r}_{\theta}(Y)$, in the same way.
Actually, the condition ``\emph{oriented}" is not necessary  for the definitions of pullbacks.
By Lemma \ref{lem fun} (3) and \cite[II. 8.1]{Br}, the pullback  defined here is compatible with the ones defined on cohomologies of  sheaves, namely,
\begin{equation}\label{compatblility1}
\rho(f^*\sigma)=f^*\rho(\sigma)
\end{equation}
for any $\sigma\in H^*_{\theta}(Y)$ (resp. $H^*_{\theta,c}(Y)$ if $f$ is proper).

Let $j:U\rightarrow X$ be the inclusion of an open subset $U$ into a (\emph{not necessarily orientable}) smooth manifold $X$.
For a sheaf  $\mathcal{F}$ on $U$,  $j_!\mathcal{F}$ is the sheaf on $X$ defined as
\begin{displaymath}
\Gamma(V,j_!\mathcal{F})=\{\,s\in\Gamma(U\cap V,\mathcal{F})\,|\,\mbox{the support of $s$ is closed relative to $V$}\,\}
\end{displaymath}
for every open subset $V$ in $X$.
By \cite[p. 184, Corollary 7.3]{I}, there is a canonical isomorphism $H^*_c(X,j_!\mathcal{F})\tilde{\rightarrow} H^*_c(U,\mathcal{F})$.
For any sheaf $\mathcal{G}$ on $X$, the adjunction morphism $j_!j^{-1}\mathcal{G}\rightarrow \mathcal{G}$
induces a morphism $H^*_c(X,j_!j^{-1}\mathcal{G})\rightarrow H^*_c(X,\mathcal{G})$.
Hence we obtain the morphism $j_!:H^*_c(U,j^{-1}\mathcal{G})\rightarrow H^*_c(X,\mathcal{G})$;
see \cite[II. 6, III. 7]{I} for more details.
In particular, we have the morphism $j_!:H^*_c(U,\underline{\mathbb{R}}_{U,\theta})\rightarrow H^*_c(X,\underline{\mathbb{R}}_{X,\theta})$.
Denote by $j_*: \mathcal{A}_c^*(U)\rightarrow\mathcal{A}_c^*(X)$ the extension  by zero, which induces the morphism $j_*: H_{\theta,c}^*(U)\rightarrow H_{\theta,c}^*(X)$.
If $X$ is \emph{oriented}, it coincides with the above pushforward $j_*$ defined by currents.

Whenever $X$ is \emph{orientable or not}, we have
\begin{prop}\label{cpb2}
Via $\rho$, $j_*$  is compatible with $j_!$ on the cohomology of sheaves, i.e., the  diagram
\begin{displaymath}
\xymatrix{
 H_{\theta,c}^*(U) \ar[d]^{j_*} \ar[r]^{\rho\quad}& H_c^*(U,\underline{\mathbb{R}}_{U,\theta})\ar[d]^{j_!}\\
 H_{\theta,c}^*(X)      \ar[r]^{\rho\quad}& H_c^*(X,\underline{\mathbb{R}}_{X,\theta}).  }
\end{displaymath}
is commutative.
\end{prop}
\begin{proof}
Denote by $\mathcal{A}_X^\bullet$ the complex $(\mathcal{A}_X^\bullet,d_\theta)$.
Let $\mathcal{A}_X^\bullet\rightarrow \mathcal{I}^\bullet$ and $j^{-1}\mathcal{I}^\bullet\rightarrow \mathcal{J}^\bullet$ be injective resolutions of complexes $\mathcal{A}_X^\bullet$ and $j^{-1}\mathcal{I}^\bullet$ of sheaves, respectively.
Then $\mathcal{I}^\bullet$ and $\mathcal{J}^\bullet$ are injective resolutions of $\underline{\mathbb{R}}_{X,\theta}$ and $\underline{\mathbb{R}}_{U,\theta}$, respectively.
Since $j_!$ is an exact functor 
, $j_!j^{-1}\mathcal{I}^\bullet\rightarrow j_!\mathcal{J}^\bullet$ is quasi-isomorphic.
By \cite[p. 41, 6.2]{I},  there exists a morphism $j_!\mathcal{J}^\bullet\rightarrow\mathcal{I}^\bullet$ of complexes such that the right triangle in the diagram
\begin{equation}\label{sheaf-comm}
\xymatrix{
 j_!\mathcal{A}_U^\bullet=j_!j^{-1}\mathcal{A}_X^\bullet\ar[d] \ar[r]& j_!j^{-1}\mathcal{I}^\bullet \ar[d]\ar[r]&j_!\mathcal{J}^\bullet\ar[dl]\\
 \mathcal{A}_X^\bullet \ar[r]& \mathcal{I}^\bullet}
\end{equation}
is commutative up to a homotopy.
Clearly, the left square in (\ref{sheaf-comm}) is commutative.
Considering  cohomologies with compact support for (\ref{sheaf-comm}), we conclude the proof.
\end{proof}

Let $j:U\rightarrow X$ be the inclusion of an open subset $U$ into an $n$-dimensional \emph{oriented} smooth manifold $X$.
For a current $T\in \mathcal{D}^{\prime p}(X)$, the current $j^*T$ is defined as $\langle j^*T,\beta\rangle=\langle T,j_*\beta\rangle$ for any $\beta\in\mathcal{A}_c^{n-p}(U)$, where $\langle ,\rangle$ is the pair of the topological dual between currents and smooth forms with compact supports.
Let $\theta$ be a closed one-form on $X$.
Clearly, $d_{\theta|_U}(j^*T)=j^*(d_\theta T)$, so $j^*$ induces $H^p_{\theta}(X)\rightarrow H^p_{\theta}(U)$.
It coincides with the pullback $j^*$ defined via forms as above,
since $\int_Uj^*\alpha\wedge\beta=\int_X\alpha\wedge j_*\beta$ for $\alpha\in\mathcal{A}^p(X)$ and $\beta\in\mathcal{A}_c^{n-p}(U)$.

\subsection{Cup product}
Let $X$ be a smooth manifold and $\theta$, $\mu$ closed  one-forms on $X$. Then
\begin{displaymath}
d_{\theta+\mu}(\beta\wedge\gamma)=d_{\theta}\beta\wedge\gamma+(-1)^{\textrm{deg}\beta}\beta\wedge d_{\mu}\gamma,
\end{displaymath}
where $\beta$ and $\gamma$ are in $\mathcal{A}^*(X)$ or $\mathcal{D}^{\prime*}(X)$, but not both in $\mathcal{D}^{\prime*}(X)$.
So we can define a \emph{cup product}
\begin{displaymath}
\cup: H^p_{\theta}(X)\times H^q_{\mu}(X)\rightarrow H^{p+q}_{\theta+\mu}(X)
\end{displaymath}
as $([\alpha]_{\theta},[\beta]_{\mu})\mapsto [\alpha\wedge\beta]_{\theta+\mu}$ for any $d_\theta$-closed $\alpha\in\mathcal{A}^p(X)$ and $d_\mu$-closed $\beta\in\mathcal{A}^p(X)$.
It can also be defined by the wedge product between  smooth forms and currents.
The two definitions coincide.
Similarly, we can define the cup products between $H^p_{\theta}(X)$ and $H^q_{\mu,c}(X)$   or $H^p_{\theta,c}(X)$ and  $H^q_{\mu,c}(X)$.
By Lemma \ref{lem fun} (2) and \cite[II. 7.5]{Br},  the cup product defined here is compatible with the one defined on  cohomologies  from sheaf theory, that is to say,
\begin{equation}\label{compatblility2}
\rho(\sigma\cup \tau)=\rho(\sigma)\cup \rho(\tau)
\end{equation}
for any $\sigma\in H^p_{\theta}(X)$ or $H^p_{\theta,c}(X)$ and $\tau\in H^q_{\mu}(X)$ or $H^q_{\mu,c}(X)$.

By \cite[Corollary 3.3.12]{Dim} and Lemma \ref{lem fun} (2),  we get the \emph{Poincar\'{e} duality theorem} for  Morse-Novikov cohomologies as follows.
\begin{cor}[{\cite[Corollary 1.4]{HR}}]\label {PD}
Let $X$ be an oriented smooth manifold of dimension $n$ and  $\theta$ a closed one-form on $X$. Then
\begin{displaymath}
PD:H_\theta^p(X)\rightarrow (H_{-\theta,c}^{n-p}(X))^*
\end{displaymath}
is an isomorphism for any $p$, where $PD([\alpha]_\theta)([\beta]_{-\theta,c})=\int_X\alpha\wedge\beta$ and $*$ denote the algebraic dual of a vector space.
\end{cor}

Let $f:X\rightarrow Y$ be a proper smooth map between oriented smooth manifolds and let $\theta$, $\mu$ be closed one-forms on $Y$.
Set $\tilde{\theta}=f^*\theta$.
Since
\begin{equation}\label{pro-f0}
f_*(T\wedge f^*\beta)=f_*T\wedge \beta
\end{equation}
for any $T\in\mathcal{D}^{\prime*}(X)$ and $\beta\in\mathcal{A}^*(Y)$,
we have the \emph{projection formula}
\begin{equation}\label{pro-f}
f_*(\sigma\cup f^*\tau)=f_*(\sigma)\cup \tau
\end{equation}
for $\sigma$ $\in$ $H_{\tilde{\theta}}^*(X)$ or $H_{\tilde{\theta},c}^*(X)$ and $\tau\in H_{\mu}^*(Y)$ or $H_{\mu,c}^*(Y)$.
If $\sigma\in H_{\tilde{\theta},c}^*(X)$ and $\tau\in H_{\mu}^*(Y)$,  the condition ``\emph{proper}" on $f$ can be removed in  $($\ref{pro-f}$)$.
By the projection formula, we get

\begin{cor}\label{i-s}
Let $f:X\rightarrow Y$ be a proper surjective smooth map of oriented  smooth manifolds with the same dimensions and $\emph{deg}f\neq 0$.
Suppose that $\theta$ is a closed one-form on $Y$ and set $\tilde{\theta}=f^*\theta$.
Then $f^*:H^*_{\theta}(Y)\rightarrow H^*_{\tilde{\theta}}(X)$ is injective and $f_*:H^*_{\tilde{\theta}}(X)\rightarrow H^{*}_{\theta}(Y)$ is surjective.
They also hold for the Morse-Novikov cohomology with compact support.
\end{cor}

\subsection{Cartesian product}
Suppose that $\theta$ and $\mu$ are closed one-forms on smooth manifolds $X$ and $Y$ respectively.
Let $pr_1$ and $pr_2$ be projections from $X\times Y$ onto $X$ and $Y$ respectively.
Set $\omega=pr_1^*\theta+pr_2^*\mu$.
The map $(\alpha,\beta)\mapsto pr_1^*(\alpha)\wedge pr_2^*(\beta)$ induces the  \emph{cartesian products}
\begin{displaymath}
\times: H^p_{\theta}(X)\times H^q_{\mu}(Y)\rightarrow H^{p+q}_{\omega}(X\times Y),
\end{displaymath}
\begin{displaymath}
\times: H^p_{\theta,c}(X)\times H^q_{\mu,c}(Y)\rightarrow H^{p+q}_{\omega,c}(X\times Y).
\end{displaymath}
By (\ref{compatblility1}) and (\ref{compatblility2}),
the cartesian products defined here are compatible with the ones defined on  cohomologies of sheaves, i.e., $\rho(\sigma\times\tau)=\rho(\sigma)\times\rho(\tau)$ for any $\sigma\in H^p_{\theta}(X)$ (resp. $H^p_{\theta,c}(X)$) and $\tau\in H^q_{\mu}(Y)$ (resp. $H^q_{\mu,c}(Y)$).

\section{K\"{u}nneth theorems}
Recall some constructions in \cite[Lemma 1.1]{HR}.
Consider the trivial bundle $\pi:\mathbb{R}^n\times F\rightarrow \mathbb{R}^n$ over $\mathbb{R}^n$, where $F$ is a smooth manifold.
Suppose that $\Theta$ is a closed one-form on $\mathbb{R}^n\times F$.
Let $pr_2:\mathbb{R}^n\times F\rightarrow F$ be the second projection and $i_0:F\rightarrow \mathbb{R}^n\times F$ an inclusion which maps $f$ to $(0,f)$.
Assume that $t$ is the coordinate of the first factor of $\mathbb{R}\times\mathbb{R}^n\times F$.
For any $p$,  define the \emph{contraction} operator $i(\partial/\partial t): \mathcal{A}^p(\mathbb{R}\times\mathbb{R}^n\times F)\rightarrow  \mathcal{A}^{p-1}(\mathbb{R}\times\mathbb{R}^n\times F)$ as
\begin{displaymath}
i(\partial/\partial t)(\Upsilon)(X_1,\ldots,X_{p-1})=\Upsilon(\partial/\partial t,X_1,\ldots,X_{p-1})
\end{displaymath}
for any $\Upsilon\in \mathcal{A}^p(\mathbb{R}\times\mathbb{R}^n\times F)$ and arbitrary smooth tangent vector fields $X_1$, $\ldots$, $X_{p-1}$ on $\mathbb{R}\times\mathbb{R}^n\times F$.
Let $g:\mathbb{R}\times\mathbb{R}^n\times F\rightarrow\mathbb{R}^n\times F$ map $(t, x, f)$ to $((1-t)x, f)$, which gives a smooth homotopy between $\textrm{id}_{\mathbb{R}^n\times F}$ and $i\circ pr_2$.
Set $u_s=\int_0^s i(\partial/\partial t)(g^*\Theta)\textrm{d}t$.
Define $K: \mathcal{A}^p(\mathbb{R}^n\times F)\rightarrow  \mathcal{A}^{p-1}(\mathbb{R}^n\times F)$ as
\begin{equation}\label{operator}
K(\alpha)=\int_0^1 e^{u_t}\cdot i(\partial/\partial t)(g^*\alpha) \textrm{d}t.
\end{equation}
Then
\begin{equation}\label{htp1}
pr_2^*i_0^*\Theta-\Theta=\textrm{d}u_1
\end{equation}
and
\begin{equation}\label{htp2}
e^{u_1}pr_2^*i_0^*\alpha-\alpha=\textrm{d}_\Theta K(\alpha)+K(\textrm{d}_\Theta \alpha)
\end{equation}
for any $\alpha\in \mathcal{A}^*(\mathbb{R}^n\times F)$, see the proof of \cite[Lemma 1.1]{HR}.

\begin{lem}\label{homotopic}
Let $\theta$ be a closed one-form on $F$ and $\tilde{\theta}=pr_2^*\theta$. Then $pr_2^*:H_\theta^*(F)\rightarrow H_{\tilde{\theta}}^*(\mathbb{R}^n\times F)$ is an isomorphism and $i_0^*$ is its inverse isomorphism.
\end{lem}
\begin{proof}
Assume that $g$, $i(\partial/\partial t)$ and $u_s$ are defined as above, where $\Theta=\tilde{\theta}$. Clearly, $i(\partial/\partial t)(g^*\tilde{\theta})=0$, and then $u_s=0$. By (\ref{htp2}), $pr_2^*\circ i_0^*=id$ on $H_{\tilde{\theta}}^*(\mathbb{R}^n\times F)$. Notice that $pr_2\circ i_0=id_F$, from which our assertion follows.
\end{proof}

\subsection{A K\"{u}nneth theorem}

\begin{thm}\label{Kun1}
Let $X$, $Y$ be smooth manifolds and let $\theta$, $\mu$ be closed one-forms on $X$, $Y$ respectively. Set $\omega=pr_1^*\theta+pr_2^*\mu$, where $pr_1$, $pr_2$ are projections from $X\times Y$ onto $X$, $Y$ respectively.
Then the cartesian product gives an isomorphism of graded vector spaces
\begin{displaymath}
H_{\theta,c}^*(X)\otimes_{\mathbb{R}}H_{\mu,c}^*(Y)\rightarrow H_{\omega,c}^*(X\times Y).
\end{displaymath}
Moreover, if $H_{\theta}^*(X)$ or $H_{\mu}^*(Y)$ has finite dimension,  the cartesian product also gives an isomorphism of graded vector spaces
\begin{displaymath}
H_{\theta}^*(X)\otimes_{\mathbb{R}}H_{\mu}^*(Y)\rightarrow H_{\omega}^*(X\times Y).
\end{displaymath}
\end{thm}

\begin{proof}
By \cite[II. 15.2]{Br},  the first part holds.
Following \cite[II. Proposition 9.12]{BT}, we easily get the second part.
\end{proof}

\begin{rem}
If $X$ is compact and $H^*_\mu(Y)$ is of finite dimension, the second part of Theorem \ref{Kun1} can be immediately obtained by \cite[IV. Theorem (15.10)]{Dem}.
\end{rem}

If $H_\theta^*(X)$ is of finite dimension, define $b_p(X,\theta)=\textrm{dim}H_\theta^p(X)$ and $\chi(X,\theta)=\sum\limits_{p\geq 0}(-1)^pb_p(X,\theta)$, which are called the $\theta$-\emph{betti number} and  $\theta$-\emph{Euler-characteristic} respectively.
Let $f:X\rightarrow X$ be a smooth self-map of an \emph{oriented compact} smooth manifold $X$ and let $\theta$ be a closed one-form on $X$ satisfying $f^*\theta=\theta$.
Then $f$ induces an endomorphism $f^*:H_\theta^*(X)\rightarrow H_\theta^*(X)$.
Define the \emph{$\theta$-Lefschetz number} of $f$ as
\begin{displaymath}
L(f,\theta)=\sum_{p\geq0}(-1)^p\textrm{tr}(f^*|_{H_\theta^p(X)}),
\end{displaymath}
where $\textrm{tr}(f^*|_{H_\theta^p(X)})$ is the trace of the endomorphism $f^*|_{H_\theta^p(X)}$.

\begin{ex}
Let $f:X\rightarrow X$ be a smooth self-map of an oriented compact smooth manifold $X$ and $\theta$ a closed one-form on $X$.
The triple $(X,\theta,f)$ satisfies the assumptions in the definition of $\theta$-Lefschetz number in the following cases:

$(1)$ $(X,0,f)$ for any $f$. In this case, $L(f,0)=L(f)$ is the classical Lefschetz number.

$(2)$ $(X,\theta,\textrm{id}_X)$ for any $\theta$. In this case, $L(\textrm{id}_X,\theta)= \chi(X,\theta)$.

$(3)$ $(X,\theta,g\cdot)$, where  $X$ is a smooth manifold with an action by a  group  $G$, $\theta$ is a $G$-invariant one-form on $X$ and $g\cdot:X\rightarrow X$ is the action on $X$ by $g\in G$.
\end{ex}

The $\theta$-betti number,  $\theta$-Euler-characteristic and $\theta$-Lefschetz number are generalizations of the corresponding concepts on the de Rham cohomology. Actually, the latter two coincide with the classical ones as follows.
\begin{prop}
$L(f,\theta)=L(f)$. In particular, $\chi(X,\theta)=\chi(X)$.
\end{prop}
\begin{proof}
Set $\textrm{dim}X=n$.
Let $\{\,[\alpha_i]_\theta\,\}$ be a basis of $H_\theta^*(X)$ and  $\{\,[\beta_j]_{-\theta}\,\}$ their dual basis in $H_{-\theta}^*(X)$ under Poincar\'{e} duality, i.e., $\int_X\alpha_i\wedge\beta_j=\delta_{ij}$, where $\alpha_i$, $\beta_j$ are all of pure degrees and $\delta_{ij}$ is the Kronecker delta.
Suppose that $pr_1$, $pr_2$ are two projections from $X\times X$ onto $X$.
Let $\Delta$ be the diagonal of $X\times X$ and $\Gamma_f$ the graph of $f$ in $X\times X$.
Let $i:\Delta\rightarrow X\times X$ and $i':\Gamma_f\rightarrow X\times X$ be inclusions and let $l:X\rightarrow\Delta$ be the diagonal map and $l':X\rightarrow \Gamma_f$ defined as $x\mapsto (x,f(x))$.
Endow $\Delta$ and $\Gamma_f$ with suitable orientations such that $l$ and $l'$ are diffeomorphisms of preserving orientations.
By Theorem \ref{Kun1}, the fundamental class $[\Gamma_f]\in H_{dR}^n(X\times X)$ can be written as $\sum\limits_{i,j} c_{ij}pr_1^*[\alpha_i]_\theta\cup pr_2^*[\beta_j]_{-\theta}$ for some $c_{ij}\in\mathbb{R}$. Set $f^*[\alpha_i]_\theta=\sum\limits_ja_{ij}[\alpha_j]_\theta$.
On one hand,
\begin{displaymath}
\int_{\Gamma_f} i'^*(pr_1^*\beta_i\wedge pr_2^*\alpha_j)=\int_X l'^*i'^*pr_1^*\beta_i\wedge l'^*i'^*pr_2^*\alpha_j=\int_X \beta_i\wedge f^*\alpha_j
=(-1)^{\textrm{deg}\alpha_i\textrm{deg}\beta_i}a_{ji}.
\end{displaymath}
On the other hand,
\begin{displaymath}
\begin{aligned}
\int_{\Gamma_f} i'^*(pr_1^*\beta_i\wedge pr_2^*\alpha_j)=&\int_{X\times X}[\Gamma_f]\wedge pr_1^*\beta_i\wedge pr_2^*\alpha_j\\
= &\sum_{k,l} c_{kl}(-1)^{(\textrm{deg}\beta_i+\textrm{deg}\alpha_j)\textrm{deg}\beta_l}\cdot\int_{X\times X}pr_1^*(\alpha_k\wedge \beta_i)\wedge pr_2^*(\alpha_j\wedge \beta_l)\\
= &(-1)^{(\textrm{deg}\beta_i+\textrm{deg}\alpha_j)\textrm{deg}\beta_j}c_{ij}.
\end{aligned}
\end{displaymath}
So $c_{ij}=(-1)^{\textrm{deg}\alpha_i\textrm{deg}\beta_i+(\textrm{deg}\beta_i+\textrm{deg}\alpha_j)\textrm{deg}\beta_j}a_{ji}$. The intersection number
\begin{displaymath}
\begin{aligned}
\Gamma_f \cdot \Delta=&\int_{X\times X}[\Gamma_f]\cup[\Delta]\\
=&(-1)^n\sum_{i,j}c_{ij}\int_\Delta i^*pr_1^*\alpha_i\wedge i^*pr_2^*\beta_j\\
= &(-1)^n\sum_{i,j}c_{ij}\int_X \alpha_i\wedge \beta_j\\
=&\sum_i (-1)^{\textrm{deg}\alpha_i}a_{ii}\\
=&\sum_p (-1)^p \textrm{tr}(f^*|_{H_\theta^p(X)})\\
=&L(f,\theta).
\end{aligned}
\end{displaymath}
Therefore, $L(f,\theta)$ is independent of $\theta$.
\end{proof}
\begin{rem}
$(1)$ For a compact complex manifold $X$, G. Bande and D. Kotschick  \cite{BK} first pointed out that $\chi(X,\theta)=\chi(X)$.
In fact, we can equip $X$ with a Riemannian metric $g$ and then define an operator $d^*_\theta$ as the formal $L^2$-adjoint of $d_\theta$ with respect to $g$.
Then $\chi(X,\theta)$ is the index of the perturbed operator $d_{\theta}+d^*_{\theta}$.
Notice that,  $[0,1]\ni t\mapsto d_{t\theta}+d^*_{t\theta}$ is a continuous family of first-order elliptic operators,  whose index is independent of $t$.

$(2)$ Let $Z$ be an oriented  submanifold  of an $n$-dimensional oriented smooth manifold $X$ with codimension $r$. Assume  that $i:Z\rightarrow X$ is the inclusion. In \cite[p. 14, (2.14)]{Dem},
$\int_Zi^*(\bullet)$ on $\mathcal{A}_c^{n-r}(X)$ defines a current on $X$,  which is closed. Its class in $H_{dR}^r(X)$ is denoted by $[Z]$. In \cite[p. 51]{BT},
the \emph{Poincar\'{e} dual} $[\eta_Z]_{dR}\in H_{dR}^r(X)$ of $Z$ is defined as $\int_Zi^*\omega=\int_X\omega\wedge\eta_Z$ for any closed $\omega\in\mathcal{A}_c^{n-r}(X)$. Then $[\eta_Z]_{dR}=(-1)^{r(n-r)}[Z]$. In this article, we use the notation $[Z]\in H_{dR}^r(X)$ and call it \emph{the fundamental class} of $Z$.
\end{rem}

\subsection{A second K\"{u}nneth theorem}
Before giving another  K\"{u}nneth theorem, we prove a lemma, which will be frequently used in what follows.
\begin{lem}\label{ind}
Let $X$ be a smooth manifold and $\mathcal{P}(U)$ a statement for any open subset $U\subseteq X$. Assume that $\mathcal{P}$ satisfies the  conditions:

$(i)$ \emph{local condition}\emph{:} There exists a basis $\mathfrak{U}$ of the topology of $X$ such that $\mathcal{P}(\bigcap\limits_{i=1}^lU_i)$ holds for any finitely many $U_1$, $\ldots$, $U_l\in \mathfrak{U}$.

$(ii)$ \emph{disjoint condition}\emph{:} Let $\{\,U_n\,|\,n\in\mathbb{N}^+\,\}$ be a collection of disjoint open subsets of $X$. If $\mathcal{P}(U_n)$ hold for all $n\in\mathbb{N}^+$, $\mathcal{P}(\bigcup\limits_{n=1}^\infty U_n)$ holds.

$(iii)$ \emph{Mayer-Vietoris condition}\emph{:} For open subsets $U$, $V$ of $X$, if $\mathcal{P}(U)$, $\mathcal{P}(V)$ and $\mathcal{P}(U\cap V)$ hold, then $\mathcal{P}(U\cup V)$ holds.\\
Then $\mathcal{P}(X)$ holds.
\end{lem}
\begin{proof}
We first prove:

$(\textbf{*})$ For open subsets  $U_1$, \ldots, $U_r$ of $X$, if $\mathcal{P}(\bigcap\limits_{j=1}^kU_{i_j})$ holds for any $1\leq i_1<\ldots<i_k\leq r$, then $\mathcal{P}(\bigcup\limits_{i=1}^r U_i)$ holds.

Obviously, (\textbf{*}) holds for $r=1$. Suppose $(\textbf{*})$ holds for $r$.
For $r+1$, set $U'_1=U_1$, \ldots, $U'_{r-1}=U_{r-1}$, $U'_r=U_r\cup U_{r+1}$.
Then $\mathcal{P}(\bigcap\limits_{j=1}^kU'_{i_j})=\mathcal{P}(\bigcap\limits_{j=1}^kU_{i_j})$ holds for any $1\leq i_1<\ldots<i_k\leq r$, where $i_k\neq r$.
Notice that $\mathcal{P}(\bigcap\limits_{j=1}^{k-1}U_{i_j}\cap U_r)$, $\mathcal{P}(\bigcap\limits_{j=1}^{k-1}U_{i_j}\cap U_{r+1})$ and  $\mathcal{P}(\bigcap\limits_{j=1}^{k-1}U_{i_j}\cap U_r\cap U_{r+1})$  hold, so does
\begin{displaymath}
\mathcal{P}(\bigcap\limits_{j=1}^{k}U'_{i_j})=\mathcal{P}((\bigcap\limits_{j=1}^{k-1}U_{i_j}\cap U_r)\cup(\bigcap\limits_{j=1}^{k-1}U_{i_j}\cap U_{r+1}))
\end{displaymath}
for any $1\leq i_1<\ldots<i_{k-1}\leq i_k=r$ by the Mayer-Vietoris condition. By the inductive hypothesis, $\mathcal{P}(\bigcup\limits_{i=1}^{r+1} U_i)=\mathcal{P}(\bigcup\limits_{i=1}^{r} U'_i)$ holds. We proved $(\textbf{*})$.

Let $\mathfrak{U}_\mathfrak{f}$ be the collection of open sets  which is the finite union of open sets in  $\mathfrak{U}$. We claim that

$(\textbf{**})$ $\mathcal{P}(V)$ holds for any finite intersection $V$ of open sets in $\mathfrak{U}_\mathfrak{f}$.

Suppose $V=\bigcap\limits_{i=1}^s U_i$, where $U_i=\bigcup\limits_{j=1}^{r_i}U_{ij}$ and $U_{ij}\in\mathfrak{U}$.
Set $\Lambda=\{\,J=(j_1,\ldots,j_s)\,|\,1\leq j_1\leq r_1,\mbox{ }\ldots, \mbox{ }1\leq j_s\leq r_s\,\}$ and $U_J=U_{1j_1}\cap\ldots\cap U_{sj_s}$. Then $V=\bigcup\limits_{J\in\Lambda}U_J$. For any $J_1$, $\ldots$, $J_t\in\Lambda$, $\mathcal{P}(U_{J_1}\cap \ldots\cap U_{J_t})$ holds by the local condition. Hence $\mathcal{P}(V)=\mathcal{P}(\bigcup\limits_{J\in\Lambda}U_J)$ holds by $(\textbf{*})$. We conclude $(**)$

By \cite[p. 16, Proposition II]{GHV},  $X=\bigcup\limits_{i=1}^lV_i$, where $V_i$ is a countable disjoint union of open sets in  $\mathfrak{U}_\mathfrak{f}$ for $1\leq i\leq l$.
For any $1\leq i_1<\ldots<i_k\leq l$, $\bigcap\limits_{j=1}^kV_{i_j}$ is a countable disjoint union of finite intersections of open sets in $\mathfrak{U}_\mathfrak{f}$.
By $(\textbf{**})$   and the disjoint condition, $\mathcal{P}(\bigcap\limits_{j=1}^kV_{i_j})$ holds, so does
$\mathcal{P}(X)$ by $(\textbf{*})$. We complete the proof.
\end{proof}

Let $\pi:E\rightarrow X$ be a smooth fiber bundle on a smooth manifold $X$. Set
\begin{displaymath}
cv=\{\,Z\subseteq E\,|\,Z\mbox{ is closed in }E\mbox{ satisfying that }\pi|_Z:Z\rightarrow X\mbox{ is proper}\,\}.
\end{displaymath}
The set $Z\in cv$ is called \emph{a compact vertical support}.  Evidently, $Z\in cv$, if and only if, $\pi^{-1}(K)\cap Z$ is compact for any compact subset $K\subseteq X$. By \cite[IV. 5.3 (b), 5.5]{Br},  $cv$ is a paracompactifying family of supports on $E$. The sheaf $\mathcal{C}_E^\infty$ of germs of smooth functions on $E$ is $cv$-soft (\cite[II. 9.4]{Br}), so is $\mathcal{A}_E^p$ for any $p$ (\cite[II. 9.16]{Br}). By \cite[II. 9.11]{Br}, $\mathcal{A}_E^p$ is $cv$-acyclic. Denote $\mathcal{A}_{cv}^*(E)=\Gamma_{cv}(E,\mathcal{A}_E^*)$. Clearly, $\mathcal{A}_{c}^*(E)\subseteq\mathcal{A}_{cv}^*(E)\subseteq\mathcal{A}^*(E)$. If $X$ is compact, $\mathcal{A}_{c}^*(E)=\mathcal{A}_{cv}^*(E)$ and if the fiber of $E$ is compact, $\mathcal{A}_{cv}^*(E)=\mathcal{A}^*(E)$.

Suppose that $\Theta$ is a closed one-form  on $E$.
By \cite[II. 4.1]{Br},  $H_{\Theta,cv}^*(E):= H^*(\mathcal{A}_{cv}^\bullet(E),\textrm{d}_{\Theta})$ can be viewed as the cohomology of $\underline{\mathbb{R}}_{E,\Theta}$ with supports in $cv$.
We call $H_{\Theta,cv}^*(E)$ the \emph{compact vertical Morse-Novikov cohomology} of $E$.
For a $d_\Theta$-closed $\alpha\in\mathcal{A}^{*}(E)$, denote by $[\alpha]_{\Theta,cv}$ its class in $H_{\Theta,cv}^*(E)$.

For any open set $U\subseteq X$, set $E_U=\pi^{-1}(U)$.
For  open subsets $U,\mbox{ }V\subseteq X$, we easily check that
\begin{equation}\label{M-V0}
\begin{aligned}
\xymatrix{
0\ar[r] &\mathcal{A}_{cv}^\bullet(E_{U\cup V})\ar[r]^{P\qquad} &\mathcal{A}_{cv}^\bullet(E_U)\oplus \mathcal{A}_{cv}^\bullet(E_V)\ar[r]^{\qquad Q} &\mathcal{A}_{cv}^\bullet(E_{U\cap V})\ar[r]&0
}
\end{aligned}
\end{equation}
is an exact sequence of complexes (following \cite[Proposition 2.3]{BT}),
where all the differentials in complexes are $d_\Theta$ and $P(\alpha)=(\alpha|_{E_U}, \alpha|_{E_V})$, $Q(\beta, \gamma)=\beta|_{E_{U\cap V}}-\gamma|_{E_{U\cap V}}$.

Denote by $E_x$ the fiber of $E$ over $x\in X$ and by $i_x:E_x\rightarrow E$ the inclusion.
For $\omega\in \mathcal{A}_{cv}^*(E)$, $\textrm{supp}(i_x^*\omega)\subseteq E_x\cap \textrm{supp}\omega$ is compact, i.e., $i_x^*\omega\in\mathcal{A}_c^*(E_x)$.  So $i_x$ induces the pullback $H_{\Theta,cv}^*(E)\rightarrow H_{\Theta|_{E_x},c}^*(E_x)$ for any closed one-form $\Theta$ on $E$.

If $E$ is an \emph{oriented} manifold, $\mathcal{A}_E^\bullet\hookrightarrow\mathcal{D}_E^{\prime\bullet}$ induces an isomorphism $H_{\Theta,cv}^*(E)\tilde{\rightarrow} H^*(\mathcal{D}_{cv}^{\prime\bullet}(E),\textrm{d}_{\Theta})$, where $\mathcal{D}_{cv}^{\prime*}(E)=\Gamma_{cv}(E,\mathcal{D}_E^{\prime*})$.
For a $d_\Theta$-closed $T\in\mathcal{D}_{cv}^{\prime*}(E)$, denote by $[T]_{\Theta,cv}$  its class in  $H_{\Theta,cv}^*(E)$.
Moreover, assume that $X$ and $E$ are both \emph{oriented} manifolds. Let $i:X\rightarrow E$ be the inclusion and $r=\textrm{rank}E$. For $T\in\mathcal{D}^{\prime*}(X)$, $i_*T\in \mathcal{D}_{cv}^{\prime*+r}(E)$. So $i_*$ induce a morphism $i_*:H_{\Theta|_X}^*(X)\rightarrow H_{\Theta,cv}^{*+r}(E)$.

Let  $\theta$, $\mu$ be closed one-forms on smooth manifolds  $X$, $Y$ respectively
and let $pr_1$, $pr_2$ be projections from $X\times Y$ onto $X$, $Y$ respectively.
The wedge product induces a cartesian product
\begin{displaymath}
H_{\theta}^*(X)\otimes_{\mathbb{R}}H_{\mu,c}^*(Y)\rightarrow H_{pr_1^*\theta+pr_2^*\mu,cv}^*(X\times Y)
\end{displaymath}
where  $X\times Y$ is viewed as a trivial fiber bundle over $X$.

\begin{lem}\label{fib-cv1}
Let $X$ be a smooth manifold and let $\theta$, $\mu$ be closed one-forms on $X$, $\mathbb{R}^n$ respectively. Set $\omega=pr_1^*\theta+pr_2^*\mu$, where $pr_1$, $pr_2$ are  projections from $X\times \mathbb{R}^n$ onto $X$, $\mathbb{R}^n$  respectively. The cartesian product
gives an isomorphism of graded vector spaces
\begin{displaymath}
H_{\theta}^*(X)\otimes_{\mathbb{R}}H_{\mu,c}^*(\mathbb{R}^n)\rightarrow H_{\omega,cv}^*(X\times \mathbb{R}^n),
\end{displaymath}
where $pr_1:X\times \mathbb{R}^n\rightarrow X$ is viewed as a smooth fiber bundle.
\end{lem}
\begin{proof}
For any open subset $U$ in $X$, denote the cartesian product by
\begin{displaymath}
\Psi^p_U:(H_{\theta}^*(U)\otimes_{\mathbb{R}}H_{\mu,c}^*(\mathbb{R}^n))^p\rightarrow H_{\omega,cv}^p(U\times \mathbb{R}^n)
\end{displaymath}
for any $p$. Denote by $\mathcal{P}(U)$ the statement that $\Psi^p_U$ is an isomorphism for any $p$.  Our goal is to prove that $\mathcal{P}(X)$ holds.  One only needs to check the three conditions in Lemma \ref{ind}. Clearly, $\mathcal{P}$ satisfies the disjoint condition.

We claim that, $\mathcal{P}(U)$ holds if $U$ is an open subset in $X$ such that $\theta|_U$ is exact.
Set $\theta|_U=dg$ for a smooth function $g$ on $U$. Since $\mathbb{R}^n$ is contractible,  $\mu=df$ for a smooth function $f$ on $\mathbb{R}^n$. We have the commutative diagram
\begin{displaymath}
\xymatrix{
 H_{\theta}^*(U)\otimes_{\mathbb{R}}H_{\mu,c}^*(\mathbb{R}^n) \ar[d]^{\cdot e^g\otimes \cdot e^f} \ar[r]^{\quad\Psi^*_U}& H_{\omega,cv}^*(U\times\mathbb{R}^n) \ar[d]^{\cdot e^{pr_1^*f+pr_2^*g}}\\
 H^*(U)\otimes_{\mathbb{R}}H_{c}^*(\mathbb{R}^n) \ar[r]^{\quad\times}&  H_{cv}^*(U\times\mathbb{R}^n), }
\end{displaymath}
where the two vertical maps are isomorphisms.
By \cite[Proposition 6.18]{BT},  the pullback of a generator of $H_c^n(\mathbb{R}^n)=\mathbb{R}$ by $pr_2^*$ is the Thom class of the vector bundle $U\times \mathbb{R}^n$ over $U$.
By \cite[Theorems 6.17, Remark 6.17.1]{BT},  the bottom row is an isomorphism.
Therefore, $\Psi^*_U$ is isomorphic, i.e., $\mathcal{P}(U)$ holds.
The claim follows.
Let $\mathfrak{U}$ be a basis of the topology of $X$ satisfying that $\theta|_U$ is exact for any $U\in\mathfrak{U}$. 
Then $\mathcal{P}(\bigcap\limits_{i=1}^lU_i)$ holds for $U_1$, \ldots, $U_l\in\mathfrak{U}$, since $\theta$ is exact on $\bigcap\limits_{i=1}^lU_i$. So $\mathcal{P}$ satisfies the local condition.

For open subsets $U$ and $V$ in $X$, there is a commutative diagram of Mayer-Vietoris sequences
\begin{displaymath}
\tiny{
\xymatrix{
   \left(H_{\theta}^*(U\cap V)\otimes_{\mathbb{R}}H_{\mu,c}^*(\mathbb{R}^n)\right)^{p-1}\ar[d]^{\Psi^{p-1}_{U\cap V}} \ar[r] &\left(H_{\theta}^*(U\cup V)\otimes_{\mathbb{R}}H_{\mu,c}^*(\mathbb{R}^n)\right)^p\ar[d]^{\Psi^p_{U\cup V}}\ar[r]&\left(H_{\theta}^*(U)\otimes_{\mathbb{R}}H_{\mu,c}^*(\mathbb{R}^n)\right)^p\oplus \left(H_{\theta}^*(V)\otimes_{\mathbb{R}}H_{\mu,c}^*(\mathbb{R}^n)\right)^p\ar[d]^{(\Psi^p_{U},\Psi^p_{V})}\cdots\\
 H_{\omega,cv}^{p-1}((U\cap V)\times \mathbb{R}^n)   \ar[r] & H_{\omega,cv}^p((U\cup V)\times \mathbb{R}^n)\ar[r]&H_{\omega,cv}^p(U\times \mathbb{R}^n)\oplus H_{\omega,cv}^p(V\times \mathbb{R}^n)\cdots,}}
\end{displaymath}
where the bottom exact sequence is induced by (\ref{M-V0}).
It implies that  $\mathcal{P}$ satisfies the Mayer-Vietoris condition by the five-lemma.
\end{proof}

View  $pr_1:X\times Y\rightarrow X$ as a trivial smooth fiber bundle and let $\Theta$ be a closed one-form on $X\times Y$.
For an open subset $U\subseteq Y$, extendings by zero give a morphism $\mathcal{A}^*_{cv}(X\times U)\rightarrow \mathcal{A}^*_{cv}(X\times Y)$.
For convenience, the image in $\mathcal{A}^*_{cv}(X\times Y)$ of $\alpha\in \mathcal{A}^*_{cv}(X\times U)$ is also denoted by $\alpha$.
If $U$ and $V$  are open subsets in $Y$, there is an exact sequence of \emph{Mayer-Vietoris type} of complexes
\begin{equation}\label{M-V}
\begin{aligned}
\small{
\xymatrix{
0\ar[r] &\mathcal{A}_{cv}^\bullet(X\times(U\cap V))\ar[r]^{P\qquad} &\mathcal{A}_{cv}^\bullet(X\times U)\oplus \mathcal{A}_{cv}^\bullet(X\times V)\ar[r]^{\qquad Q} &\mathcal{A}_{cv}^\bullet(X\times (U\cup V))\ar[r]&0
},}
\end{aligned}
\end{equation}
where all the differentials in complexes are $d_\Theta$ and $P(\alpha)=(\alpha, -\alpha)$, $Q(\beta, \gamma)=\beta+\gamma$.
One can check it as the one for $\mathcal{A}_c^\bullet(\bullet)$, refer to \cite[Proposition 2.7]{BT}.

Assume that $Y=\bigsqcup\limits_{\alpha\in I} Y_\alpha$ is a disjoint union of smooth manifolds.
For  a form $\omega$ on $X\times Y$ and any compact subset $K\subseteq X$, $(K\times Y)\cap \textrm{supp}\omega$ is compact if and only if $(K\times Y_\alpha)\cap \textrm{supp}\omega$ is nonempty for only finitely many $\alpha\in I$ and they are all compact. So
\begin{displaymath}
\mathcal{A}_{cv}^*(X\times Y)=\bigoplus_{\alpha\in I}\mathcal{A}_{cv}^*(X\times Y_\alpha),
\end{displaymath}
where $X\times Y$ and $X\times Y_\alpha$ are viewed as smooth fiber bundles over $X$.

Now, we give another K\"{u}nneth theorem as follows.
\begin{thm}\label{Kun2}
Let $X$ and $Y$ be smooth manifolds and let $\theta$, $\mu$ be closed one-forms on $X$, $Y$ respectively. Set $\omega=pr_1^*\theta+pr_2^*\mu$, where $pr_1$, $pr_2$ are projections from $X\times Y$ onto $X$, $Y$ respectively. The cartesian product gives an isomorphism of graded vector spaces
\begin{displaymath}
H_{\theta}^*(X)\otimes_{\mathbb{R}}H_{\mu,c}^*(Y)\rightarrow H_{\omega,cv}^*(X\times Y),
\end{displaymath}
where $pr_1:X\times Y\rightarrow X$ is viewed as a smooth fiber bundle.
\end{thm}
\begin{proof}
For any open set $U$ in $Y$, the cartesian product is denoted by
\begin{displaymath}
\Psi^p_U:(H_{\theta}^*(X)\otimes_{\mathbb{R}}H_{\mu,c}^*(U))^p\rightarrow H_{\omega,cv}^p(X\times U)
\end{displaymath}
for any $p$. Let $\mathcal{P}(U)$ be the statement that $\Psi^p_U$ is an isomorphism for any $p$.  The theorem is equivalent to saying that $\mathcal{P}(Y)$ holds.  We only need to check the three conditions in Lemma \ref{ind}. Clearly, $\mathcal{P}$ satisfies the disjoint condition.
Let $\mathfrak{U}$ be a basis of the topology of $Y$ such that it is a good covering of $Y$.
For any $U_1$, \ldots, $U_l\in\mathfrak{U}$, $\bigcap\limits_{i=1}^lU_i$ is diffeomorphic to $\mathbb{R}^m$, where $m=\textrm{dim}Y$.
By Lemma \ref{fib-cv1}, $\mathcal{P}(\bigcap\limits_{i=1}^lU_i)$ holds, so $\mathcal{P}$ satisfies the local condition.
By a diagram of Mayer-Vietoris sequences and the five-lemma, $\mathcal{P}$ satisfies the Mayer-Vietoris condition.
\end{proof}

\section{Leray-Hirsch theorems}
\begin{lem}\label{cv}
Let $F$ be a smooth manifold and let $K: \mathcal{A}^*(\mathbb{R}^n\times F)\rightarrow  \mathcal{A}^{*-1}(\mathbb{R}^n\times F)$ be defined as \emph{(\ref{operator})}. If $\alpha\in \mathcal{A}_{cv}^*(\mathbb{R}^n\times F)$, then $K(\alpha)\in \mathcal{A}_{cv}^{*-1}(\mathbb{R}^n\times F)$.
\end{lem}
\begin{proof}
Let $pr_{23}:\mathbb{R}\times\mathbb{R}^n\times F\rightarrow \mathbb{R}^n\times F$ be the projection map and  $g:\mathbb{R}\times\mathbb{R}^n\times F\rightarrow\mathbb{R}^n\times F$ map $(t, x, f)$ to $((1-t)x, f)$. Set
\begin{displaymath}
C=pr_{23}\left(g^{-1}(\textrm{supp}\alpha)\cap ([0,1]\times\mathbb{R}^n\times F)\right).
\end{displaymath}
Obviously, $C$ is closed in $\mathbb{R}^n\times F$ and $g^{-1}(\textrm{supp}\alpha)\cap ([0,1]\times\mathbb{R}^n\times F)\subseteq [0,1]\times C$.
By \cite[p. 179, Remark]{GHV}, $\textrm{supp}\left(i(\partial/\partial t)g^*\alpha\right)\cap ([0,1]\times\mathbb{R}^n\times F)\subseteq [0,1]\times C$.
So $\textrm{supp}K(\alpha)\subseteq C$.

For arbitrary compact subset $L\subseteq \mathbb{R}^n$, let $h:[0,1]\times L\rightarrow\mathbb{R}^n$ map $(t, x)$ to $(1-t)x$.
Then $h\times \textrm{id}_F$ is the restriction of $g$ on $[0,1]\times L\times F$. Choose a compact ball $B\supseteq L$ with the center at the original point. Then
\begin{displaymath}
\begin{aligned}
\pi^{-1}(L)\cap \textrm{supp}K(\alpha)\subseteq& (L\times F)\cap pr_{23}\left(g^{-1}(\textrm{supp}\alpha)\cap ([0,1]\times\mathbb{R}^n\times F)\right)\\
=&pr_{23}\left(([0,1]\times L\times F)\cap g^{-1}(\textrm{supp}\alpha)\right)\\
\subseteq &pr_{23}\left((h\times id_F)^{-1}\left((B\times F)\cap\textrm{supp}\alpha\right)\right).
\end{aligned}
\end{displaymath}
Clearly, $h\times id_F$ is proper and $(B\times F)\cap\textrm{supp}\alpha$ is compact,  so $\pi^{-1}(L)\cap \textrm{supp}K(\alpha)$ is compact. We complete the proof.
\end{proof}

We give a theorem of Leray-Hirsch type on Morse-Novikov cohomologies, which will be used to compute the cohomologies of projective bundles.
\begin{thm}\label{L-H}
Let $\pi:E\rightarrow X$ be a smooth fiber bundle over a smooth manifold $X$ and let $\theta$, $\Omega$ be closed one-forms on $X$, $E$ respectively. Set $\tilde{\theta}=\pi^*\theta$.

$(1)$ Assume that  there exist classes $e_1$, $\dots$, $e_r$ of pure degrees in $ H_{\Omega}^*(E)$ such that  their restrictions $e_1|_{E_x}$, $\dots$, $e_r|_{E_x}$ freely linearly generate $H_{\Omega|_{E_x}}^*(E_x)$ for every $x\in X$. Then
\begin{displaymath}
\sum\limits_{i=1}^r\pi^*(\bullet)\cup e_i:\bigoplus_{i=1}^rH_{\theta}^{*-u_i}(X)\rightarrow H_{\tilde{\theta}+\Omega}^*(E)
\end{displaymath}
is an isomorphism, where $\emph{deg} e_i=u_i$ for $1\leq i\leq r$.

$(2)$ Assume that there exist classes $e_1$, $\dots$, $e_r$ of pure degrees in $ H_{\Omega,cv}^*(E)$ such that their restrictions $e_1|_{E_x}$, $\dots$, $e_r|_{E_x}$ freely linearly generate $H_{\Omega|_{E_x},c}^*(E_x)$ for every $x\in X$. Then $\sum\limits_{i=1}^r\pi^*(\bullet)\cup e_i$ gives two isomorphisms
\begin{displaymath}
\bigoplus_{i=1}^rH_{\theta,c}^{*-u_i}(X) \tilde{\rightarrow} H_{\tilde{\theta}+\Omega,c}^*(E),
\end{displaymath}
\begin{displaymath}
\bigoplus_{i=1}^rH_\theta^{*-u_i}(X) \tilde{\rightarrow} H_{\tilde{\theta}+\Omega,cv}^*(E),
\end{displaymath}
where $\emph{deg} e_i=u_i$ for $1\leq i\leq r$.

$(3)$ Assume that there exist classes $e_1$, $\dots$, $e_r$ of pure degrees in $H_{\Omega,c}^*(E)$, such that  their restrictions $e_1|_{E_x}$, $\dots$, $e_r|_{E_x}$ freely linearly generate $H_{\Omega|_{E_x},c}^*(E_x)$ for every $x\in X$. Then
\begin{displaymath}
\sum\limits_{i=1}^r\pi^*(\bullet)\cup e_i:\bigoplus_{i=1}^rH_{\theta,c}^{*-u_i}(X) \rightarrow H_{\tilde{\theta}+\Omega,c}^*(E)
\end{displaymath}
is an isomorphism, where $\emph{deg} e_i=u_i$ for $1\leq i\leq r$.
\end{thm}
\begin{proof}
We only prove the second result in $(2)$ and the others can be obtained similarly.

Let $F$ be the general fiber of $E$ and $\textrm{dim}X=n$.
For any open set $U\subseteq X$, set
\begin{displaymath}
\Psi^p_U=\sum\limits_{i=1}^r\pi^*(\bullet)\cup e_i: \bigoplus\limits_{i=1}^rH_{\theta,c}^{p-u_i}(U)\rightarrow H_{\tilde{\theta}+\Omega,cv}^p(E_U),
\end{displaymath}
where $E_U=\pi^{-1}(U)$. Let $\mathcal{P}(U)$ be the statement that $\Psi^p_U$ is an isomorphism for any $p$.
We aim to show that $\mathcal{P}(X)$ holds.  It suffices to check the three conditions in Lemma \ref{ind}.
Clearly, $\mathcal{P}$ satisfies the disjoint condition.
By a commutative diagram of Mayer-Vietoris sequences and the five-lemma, $\mathcal{P}$  satisfies the Mayer-Vietoris condition.

We claim that, $\mathcal{P}(U)$ holds if the open subset  $U\subseteq X$  is diffeomorphic to $\mathbb{R}^n$ such that $E_U$ is smooth trivial.
Let $\varphi_U: U\times F\rightarrow E_U$ be a smooth trivialization of $E$ on $U$ and let $pr_1$, $pr_2$ be projections from $U\times F$ onto $U$, $F$ respectively, which satisfy $\pi\circ \varphi_U=pr_1$.
Fixed a point $o\in U$, set $j_o:F\rightarrow U\times F$ as $f\mapsto (o,f)$.
Clearly, $pr_2\circ j_o=id_F$ and $i_o:=\varphi_U\circ j_o$ is the embedding $F\hookrightarrow E_U$ of the fiber $E_o$ over $o$ into $E_U$.
Set $e_i^{\prime}=(\varphi_U)^*e_i$ in $H_{\varphi_U^*\Omega,cv}^*(U\times F)$, $i=1$, $\ldots$, $r$.
Let $\{\beta_i\}_{i=1}^r$ be a system of $d_{\varphi_U^*\Omega}$-closed forms of pure degrees in $\mathcal{A}_{cv}^*(U\times F)$, such that $e'_i=[\beta_i]_{\varphi_U^*\Omega,cv}$ for $1\leq i\leq r$.
Then $j_o^*e_i^{\prime}=i_o^*e_i$ for any $i$.
The linear independence of $i_o^*e_1$, \ldots, $i_o^*e_r$ implies that $e'_1$, \ldots, $e'_r$ are also linearly independent, so mapping $e_i$ to $e'_i$ for $1\leq i\leq r$ naturally give a linear isomorphism $\textrm{span}_{\mathbb{R}}\{e_1, \ldots, e_r\}\rightarrow\textrm{span}_{\mathbb{R}}\{e_1^{\prime}, \ldots, e_r^{\prime}\}$.
By (\ref{htp1}), (\ref{htp2}) and Lemma \ref{cv}, there exists a smooth function $u$ on $U\times F$ such that
\begin{equation}\label{subtle-1}
pr_2^*i_0^*\Omega-\varphi_U^*\Omega=pr_2^*j_0^*\varphi_U^*\Omega-\varphi_U^*\Omega=du,
\end{equation}
\begin{equation}\label{subtle-2}
\beta_i-e^{u}\cdot pr_2^*j_0^*\beta_i=d_{\varphi_U^*\Omega}\gamma_i
\end{equation}
for some $\gamma_i\in \mathcal{A}_{cv}^*(U\times F)$. Then
\begin{equation}\label{subtle-3}
e^{-u}\cdot pr_1^*\alpha\wedge \beta_i=pr_1^*\alpha\wedge pr_2^*j_0^*\beta_i+(-1)^{\textrm{deg}\alpha}d_{pr_1^*\theta+pr_2^*i_0^*\Omega}(e^{-u}\cdot pr_1^*\alpha\wedge\gamma_i)
\end{equation}
for any $d_\theta$-closed form $\alpha\in\mathcal{A}^*(U)$.
Then $pr_1^*\alpha\wedge\gamma_i$ has a compact support for any  $\alpha\in\mathcal{A}_{c}^*(U)$.
There is a commutative diagram
\begin{displaymath}
\xymatrix{
                & \bigoplus\limits_{i=1}^rH_{\theta}^{p-u_i}(U) \ar[dl]_{\Psi^p_U}\ar[dr]^{\quad\quad\sum\limits_{i=1}^rpr_1^*(\bullet)\cup e'_i}   \ar[rr]^{\sum\limits_{i=1}^r\bullet\otimes j_0^*e'_i} & &     (H_{\theta}^*(U)\otimes_{\mathbb{R}} H_{i_0^*\Omega,c}^*(F))^p \ar[d]^{\times}           \\
 H_{\tilde{\theta}+\Omega,cv}^p(E_U)  \ar[rr]_{\cong}^{\varphi_U^*} &&     H_{pr_1^*\theta+\varphi_U^*\Omega,cv}^p(U\times F)  \ar[r]_{\cong}^{\cdot e^{-u}} &     H_{pr_1^*\theta+pr_2^*i_0^*\Omega,cv}^p(U\times F)  }
\end{displaymath}
for any $p$, where the top map is clearly an isomorphism.
By  Theorem \ref{Kun2}, the vertical map  is an isomorphism, so is $\Psi^p_U$. The claim is verified. Let $\mathfrak{U}$ be a basis of the topology of $X$ such that it is a good covering of $X$  and $E_U$ is smooth trivial for any $U\in\mathfrak{U}$. For $U_1$, \ldots, $U_l\in\mathfrak{U}$, $\bigcap\limits_{i=1}^lU_i$ is diffeomorphic to $\mathbb{R}^n$ and $E_{U_1\cap\ldots U_l}$ is smooth trivial, so $\mathcal{P}(\bigcap\limits_{i=1}^lU_i)$ holds. Hence $\mathcal{P}$ satisfies the local condition.
\end{proof}
\begin{rem}
For $\Omega=0$,  H. Haddou \cite{Ha} proved Theorem \ref{L-H} (1)  for the case that $X$ has a finite good covering and Y. Kawahara \cite{K} gave a holomorphic version of Theorem \ref{L-H} (1) in the category of complex affine manifolds.
\end{rem}

Let $\mathbb{P}(E)$ be the complex projectivization of a complex smooth vector bundle $E$ of complex rank $r$ on a smooth manifold $X$.
Then $\mathcal{O}_{\mathbb{P}(E)}(-1)=\{(l,v)\in \mathbb{P}(E)\times E|\mbox{ }v\in l\}$ is a complex line bundle over $\mathbb{P}(E)$, which is said to be the \emph{universal line bundle} over $\mathbb{P}(E)$.
\begin{cor}\label{3}
Let $\pi:\mathbb{P}(E)\rightarrow X$ be the complex projectivization of a complex smooth vector bundle $E$ of complex rank $r$ on a smooth manifold $X$ and
let $\theta$ be a closed one-form on $X$.
Assume that $\tilde{\theta}=\pi^*\theta$ and $h=c_1(\mathcal{O}_{\mathbb{P}(E)}(-1))\in H_{dR}^2({\mathbb{P}(E)})$ is the first Chern class of $\mathcal{O}_{\mathbb{P}(E)}(-1)$.
Then $\sum\limits_{i=0}^{r-1}\pi^*(\bullet)\cup h^i$ gives two isomorphisms
\begin{equation}\label{proj-bun-1}
\bigoplus\limits_{i=0}^{r-1}H_{\theta}^{*-i}(X)\tilde{\rightarrow}H_{\tilde{\theta}}^*(\mathbb{P}(E)),
\end{equation}
\begin{equation}\label{proj-bun-2}
\bigoplus\limits_{i=0}^{r-1}H_{\theta,c}^{*-i}(X)\tilde{\rightarrow}H_{\tilde{\theta},c}^*(\mathbb{P}(E)).
\end{equation}
\end{cor}
\begin{proof}
For every $x\in X$, $1$, $h$, \ldots, $h^{r-1}$ restricted to the fibre $\pi^{-1}(x)=\mathbb{P}(E_x)$ is a basis of $H_{dR}^*(\mathbb{P}(E_x))$.  By Theorem \ref{L-H} $(1)$ $(2)$, we proved the conclusion.
\end{proof}

Suppose $\pi:E\rightarrow X$ is an \emph{oriented} smooth vector bundle of rank $r$ on a (\emph{not necessarily orientable}) smooth manifold $X$ and $\omega\in \mathcal{A}_{cv}^p(E)$.
For a chart $U$ on $X$ satisfying that  $E_U$ is trivial, let $(x_1,\ldots,x_n;t_1,\ldots,t_r)$ be the local coordinates of $E$ such that $dt_1\wedge\ldots\wedge dt_r$ gives the orientation of $E$.
If $\omega=\sum\limits_{|I|+|J|=p} f_{I,J}(t,x)dt_I\wedge dx_J$ on $E_U$, then
\begin{displaymath}
\sum_{|J|=p-r}\left(\int_{\mathbb{R}^r}f_{1\ldots r,J}(t,x)dt_1\wedge\ldots\wedge dt_r\right)dx_J
\end{displaymath}
defines a $(p-r)$-form on $U$. For various charts of $X$, these local forms give a global one on $X$, denoted by $\pi_*\omega$.

\begin{rem}
Up to the sign $(-1)^{r(p-r)}$, $\pi_*\omega$ defined here coincides with the one defined in \cite[p. 61-62]{BT}. Moreover, if $X$ is \emph{oriented},  $\pi_*\omega$ defined here  is just the pushforward of $\omega$ as currents \cite[p. 18, (2.16)]{Dem}.
\end{rem}
By Thom isomorphism theorem \cite[Theorem 12.2, (12.2.1)]{BT},  $\pi_*:H_{cv}^*(E)\rightarrow H_{dR}^{*-r}(X)$ is an isomorphism.
Let $\Phi\in\mathcal{A}_{cv}^r(E)$ satisfy $\pi_*[\Phi]_{cv}=1$ in $H_{dR}^0(X)=\mathbb{R}$.
Then $[\Phi]_{cv}\in H_{cv}^r(E)$ is  the \emph{Thom class} of $E$.
Evidently, $\pi_*\Phi=1$ in $\mathcal{A}^0(X)$.
In addition, assume that $X$ is an \emph{oriented} smooth manifold.
Let $i:X\rightarrow E$ be the inclusion of the zero section of $E$ and $r=\textrm{rank}E$.
For $T\in\mathcal{D}^{\prime*}(X)$, $i_*T\in \mathcal{D}_{cv}^{\prime*+r}(E)$. So $i_*$ induce a morphism $i_*:H_{\Theta|_X}^*(X)\rightarrow H_{\Theta,cv}^{*+r}(E)$.

\begin{cor}[Thom isomorphism]\label{is}
Let $\pi:E\rightarrow X$  be an oriented smooth vector bundle of rank $r$ on a smooth manifold $X$. Assume that $\theta$ is a closed one-form on $X$. Then $[\Phi]_{cv}\cup\pi^*(\bullet)$ gives isomorphisms
$H_{\theta,c}^{*-r}(X)\tilde{\rightarrow} H_{\tilde{\theta},c}^*(E)$
and
$H_{\theta}^{*-r}(X)\tilde{\rightarrow} H_{\tilde{\theta},cv}^*(E)$,
which have the inverse isomorphism $\pi_*$. Moreover, if $X$ is oriented, $[\Phi]_{cv}\cup\pi^*(\bullet)$ coincides with the pushforward $i_*$ on both cases, where  $i:X\rightarrow E$ is the inclusion of the zero section of $E$.
\end{cor}
\begin{proof}
By \cite[Proposition 6.18]{BT},  the restriction $[\Phi]_{dR}|_{E_x}$ is a generator of $H_{dR,c}^*(E_x)$.
By Theorem \ref{L-H} $(2)$, $[\Phi]_{cv}\cup\pi^*(\bullet)$ gives the two isomorphisms.
For every $\alpha\in\mathcal{A}^*(X)$, $\pi_*(\Phi\wedge\pi^*\alpha)=\alpha$, so
$\pi_*:H_{\tilde{\theta},c}^*(E)\rightarrow H_{\theta,c}^{*-r}(X)$  and $\pi_*:H_{\tilde{\theta},cv}^*(E)\rightarrow H_{\theta}^{*-r}(X)$ are their inverse isomorphisms.
If $X$ is oriented, $i_*$ is  well-defined and $\pi_*i_*=id$.
So $i_*=\pi_*^{-1}=[\Phi]_{cv}\cup\pi^*(\bullet)$.
\end{proof}

\section{Blow-up formulae}
Now, we  prove Lemma \ref{key} as follows.
\begin{proof}
Set $r=\textrm{codim}Y$.
Let $N\cong N_{Y/X}$ (as smooth manifolds) be a tubular neighborhood of $Y$ in $X$ and denote by $\tau$ the projection of the vector bundle $N$ over $Y$.
Let $l:Y\rightarrow N$ and $j:N\rightarrow X$ be inclusions.
Denote by $[Y]_N\in H_{dR}^r(N)$ the fundamental class of $Y$ in $N$.
Notice that $j$ is smoothly homotopic to $i\circ \tau$. So $\theta|_N-\tau^*(\theta|_Y)=du$ for some $u\in\mathcal{A}^0(N)$.
Since $\tau\circ l=id_Y$,  $d(l^*u)=0$, i.e., $l^*u=c$ is a constant.
Replace $u$ with $u-c$, then $l^*u=0$.
By (\ref{pro-f0}),
\begin{equation}\label{im}
e^{-u}\cdot l_*T=l_*(e^{-l^*u}\cdot T)=l_*T,
\end{equation}
for any current $T$ on $Y$.

By the localization principle (\cite[Proposition 6.25]{BT}), there is a representative $\eta_Y\in\mathcal{A}^r(X)$ of $[Y]$ such that $\textrm{supp}\eta_Y\subseteq N$.
Then $[\eta_Y|_N]_{dR}=[Y]_N$.
We have $l_*(1)=\eta_Y|_N+dS$ for some $S\in\mathcal{D}^{\prime r-1}(N)$.
Let $\alpha\in \mathcal{A}^*(Y)$ be a representative of $\sigma\in H_{\theta|_Y}^*(Y)$.
As a current on $N$,
\begin{displaymath}
\begin{aligned}
l_*\alpha=&e^{-u}\cdot l_*\alpha=e^{-u}\cdot l_*l^*\tau^*\alpha\\
=&e^{-u}\cdot l_*(1)\wedge\tau^*\alpha\\
=&e^{-u}\cdot \eta_Y|_N\wedge\tau^*\alpha+d_{\theta|_N}(e^{-u}\cdot S\wedge\tau^*\alpha).
\end{aligned}
\end{displaymath}
Since $\textrm{supp}(l_{*}\alpha)\subseteq Y$, $j|_{\textrm{supp}(l_{*}\alpha)}$ is proper, which implies that $j_*(l_{*}\alpha)$ is defined well.  Clearly, $j_*(l_{*}\alpha)=i_{*}\alpha$ and $j^*j_*(l_{*}\alpha)=l_{*}\alpha$. Thus
\begin{displaymath}
\begin{aligned}
i^*i_*\sigma=&l^*j^*[i_*\alpha]_\theta=l^*[j^*j_*l_*\alpha]_{\theta|_N}\\
=&l^*[l_*\alpha]_{\theta|_N}=l^*[e^{-u}\cdot \eta_Y|_N\wedge\tau^*\alpha]_{\theta|_N}\\
=&[Y]|_Y\cup\sigma.
\end{aligned}
\end{displaymath}

Let $\Phi\in\mathcal{A}_{cv}^r(N)$ be a representative of the Thom class of the vector bundle $N$ satisfying $\tau_*\Phi=1$.
Let $\alpha\in \mathcal{A}_c^*(Y)$ be a representative of $\sigma\in H_{\theta|_Y,c}^*(Y)$.
By Lemma \ref{is},
\begin{equation}\label{im2}
l_*\alpha=\Phi\wedge\tau^*\alpha+d_{\tau^*(\theta|_Y)}S
\end{equation}
for some  $S\in\mathcal{D}_c^{\prime p+r-1}(N)$.
Combining (\ref{im}) and (\ref{im2}), we have
\begin{displaymath}
l_*\alpha=e^{-u}\cdot\Phi\wedge\tau^*\alpha+d_{\theta|_N}(e^{-u}S).
\end{displaymath}
Therefore,
\begin{displaymath}
\begin{aligned}
i^*i_*\sigma=&i^*j_*[l_*\alpha]_{\theta|_N,c}=i^*j_*[e^{-u}\cdot\Phi\wedge\tau^*\alpha]_{\theta|_N,c}\\
=&[l^*j^*j_*(e^{-u}\Phi\wedge\tau^*\alpha)]_{\theta|_Y,c}\\
=&[\Phi]|_Y\cup\sigma,
\end{aligned}
\end{displaymath}
where we used that $j^*j_*=id$ on $\mathcal{A}_c^*(N)$. By \cite[I. Proposition 6.24 (b)]{BT}, $[\Phi]_{dR}=[Y]_N$ in $H_{dR}^r(N)$. Since $[Y]_N|_Y=[Y]|_Y$, $i^*i_*\sigma=[Y]|_Y\cup\sigma$.

We complete the proof.
\end{proof}

\begin{lem}\label{1}
Let $\pi:E\rightarrow X$ be a smooth vector bundle of rank $r$ on a smooth manifold $X$ and $\Theta$ a closed one-form on $E$.
Let $U\subseteq E$ be an open neighborhood of the zero section of $E$  and denote by $i_U:X\rightarrow U$ the inclusion of the zero section of $E$ into $U$.
Then $i_U^*:H_{\Theta}^*(U)\tilde{\rightarrow} H_{\Theta|_X}^*(X)$ is an isomorphism.
Moreover, if $X$ and $E$ are oriented, then $i_*: H_{\Theta|_X,c}^*(X)\tilde{\rightarrow}H_{\Theta,c}^{*+r}(E)$ and $i_*: H_{\Theta|_X}^*(X)\tilde{\rightarrow}H_{\Theta,cv}^{*+r}(E)$ are isomorphisms, where $i:X\rightarrow E$ is the inclusion of the zero section.
\end{lem}
\begin{proof}
Since $\pi|_U\circ i_{U}=id_X$, $i_{U}^*\circ(\pi|_U)^*=id:H_{\Theta|_X}^*(X)\rightarrow H_{\pi^*(\Theta|_X)}^*(U)\rightarrow H_{\Theta|_X}^*(X)$.
Define $g:U\times [0,1]\rightarrow U$ as $(e,t)\mapsto t\cdot e$.
Then $g(\cdot,0)=i_{U}\circ \pi|_U$ and $g(\cdot,1)=id_U$, i.e., $g$ gives a smooth homotopy between $i_{U}\circ \pi|_U$ and $\textrm{id}_U$.
Clearly, $i_U\circ \pi|_U\circ g (e,t)=i_U \circ\pi|_U(e)$ is independently with $t$, so $i(\partial/\partial t)\left(g^*(\pi|_U)^*(\Theta|_X)\right)=0$.
By \cite[Lemma 1.1]{HR}, $(\pi|_U)^*\circ i_{U}^*=id:H_{\pi^*(\Theta|_X)}^*(U)\rightarrow H_{\Theta|_X}^*(X)\rightarrow H_{\pi^*(\Theta|_X)}^*(U)$.
Hence $i_U^*:H_{\pi^*(\Theta|_X)}^*(U)\rightarrow H_{\Theta|_X}^*(X)$ is an isomorphism.
With similar arguments as in the proof of Lemma \ref{key}, we can choose $u\in\mathcal{A}^0(U)$ such that $\Theta|_U-\pi^*(\Theta|_X)|_U=\Theta|_U-(\pi|_U)^*i_U^*(\Theta|_U)=du$ and $i_U^*u=0$.
By (\ref{isom}), $e^u\cdot:H_{\Theta}^*(U)\rightarrow H_{\pi^*(\Theta|_X)}^*(U)$ is an isomorphism.
The diagram
\begin{displaymath}
\xymatrix{
  H_{\Theta}^*(U) \ar[rr]_{e^u\cdot}^{\cong} \ar[dr]_{i_U^*}
                &  &    H_{\pi^*(
                \Theta|_X)}^*(U) \ar[dl]_{\cong}^{i_U^*}    \\
                & H_{\Theta|_X}^*(X) }
\end{displaymath}
is commutative, which implies that $i_U^*:H_{\Theta}^p(U)\rightarrow H_{\Theta|_X}^p(X)$ is an isomorphism.

Assume that $X$ and $E$ are oriented.
By (\ref{pro-f0}), $e^u\cdot i_*\alpha=i_*(e^{i^*u}\cdot \alpha)=i_*\alpha$ for any $\alpha\in \mathcal{A}_c^*(X)$.
There is a commutative diagram
\begin{displaymath}
\xymatrix@R=0.5cm{
                &     H_{\Theta,c}^{*+r}(E)     \ar[dd]_{e^u\cdot}^{\cong}     \\
 H_{\Theta|_X,c}^*(X) \ar[ur]^{i_*} \ar[dr]_{i_*}                 \\
                &        H_{\pi^*(
                \Theta|_X),c}^{*+r}(E).                 }
\end{displaymath}
By (\ref{isom2}) and  Corollary \ref{is}, $e^u\cdot$  and $i_*: H_{\Theta|_X,c}^*(X)\rightarrow H_{\pi^*(\Theta|_X),c}^{*+r}(E)$ are isomorphisms, so is $i_*: H_{\Theta|_X,c}^*(X)\rightarrow H_{\Theta,c}^{*+r}(E)$.
In the same way, $i_*: H_{\Theta|_X}^*(X)\rightarrow H_{\Theta,cv}^{*+r}(E)$ is also an isomorphism.
\end{proof}

Now, we give a proof of Theorem \ref{1.3}.
\begin{proof}
Set $U=X-Y$ and $\widetilde{U}=\widetilde{X}-E$. Then $\pi|_{\widetilde{U}}:\widetilde{U}\rightarrow U$ is biholomorphic.
Since $\pi$ is proper, we can choose a tubular neighborhood $V\subseteq X$ of $Y$ such that $\widetilde{V}=\pi^{-1}(V)$ is contained in a tubular neighborhood of $E$ in $\widetilde{X}$.
Set $W=U\cap V$ and $\widetilde{W}=\widetilde{U}\cap\widetilde{V}$. Then $\pi|_{\widetilde{W}}:\widetilde{W}\rightarrow W$ is biholomorphic.
There is a commutative diagram of Mayer-Vietoris sequences
\begin{displaymath}
\small{\xymatrix{
    \cdots\ar[r]&H_{\theta}^{k-1}(W)\ar[d]^{\cong}\ar[r]&H_{\theta}^k(X) \ar[d]^{\pi^*} \ar[r]& H_{\theta}^k(U)\oplus H_{\theta}^k(V) \ar[d]^{(\pi|_{\widetilde{U}})^*\oplus(\pi|_{\widetilde{V}})^*}\ar[r]&  H_{\theta}^{k}(W)\ar[d]^{\cong}\ar[r]&H_{\theta}^{k+1}(X)\ar[d]^{\pi^*}\ar[r]&\cdots\\
 \cdots\ar[r]&H_{\tilde{\theta}}^{k-1}(\widetilde{W})     \ar[r] & H_{\tilde{\theta}}^k(\widetilde{X})\ar[r]&H_{\tilde{\theta}}^k(\widetilde{U})\oplus H_{\tilde{\theta}}^k(\widetilde{V})       \ar[r]& H_{\tilde{\theta}}^k(\widetilde{W})     \ar[r] & H_{\tilde{\theta}}^{k+1}(\widetilde{X})\ar[r]&\cdots     .}}
\end{displaymath}
By Corollary \ref{i-s}, $\pi^*$ is injective. By the snake lemma (\cite[p. 4]{I}), the restrictions induce an isomorphism
\begin{equation}\label{*}
\textrm{coker}\pi^*\tilde{\rightarrow}\textrm{coker}\left((\pi|_{\widetilde{U}})^*\oplus(\pi|_{\widetilde{V}})^*\right)\cong\textrm{coker}(\pi|_{\widetilde{V}})^*.
\end{equation}
Let $i'_Y:Y\rightarrow V$ and $i'_E:E\rightarrow \widetilde{V}$ be the inclusions. By Lemma \ref{1}, $i_Y^{\prime*}: H_{\theta}^k(V)\rightarrow H_{\theta|Y}^k(Y)$ and  $i_E^{\prime*}:H_{\tilde{\theta}}^k(\widetilde{V})\rightarrow H_{\tilde{\theta}|_E}^k(E)$ are isomorphisms. Since $\pi|_{\widetilde{V}}\circ i'_E=i'_Y\circ\pi|_E$, $i_E^{\prime*}$ induces an isomorphism
\begin{equation}\label{**}
\textrm{coker}(\pi|_{\widetilde{V}})^*\tilde{\rightarrow}\textrm{coker}(\pi|_E)^*.
\end{equation}
By Corollary \ref{3}, $(\pi|_E)^*$ is injective.
Denote by $i_Y:Y\rightarrow X$ the inclusion.
Combining (\ref{*}) and (\ref{**}), we have a commutative diagram of short exact sequences
\begin{equation}\label{commutative1}
\xymatrix{
 0\ar[r]&H_\theta^k(X)\ar[d]^{i_Y^*} \ar[r]^{\pi^*}& H_{\tilde{\theta}}^k(\widetilde{X})\ar[d]^{i_E^*} \ar[r]& \textrm{coker}\pi^* \ar[d]^{\cong}\ar[r]& 0\\
 0\ar[r]&H_{\theta|_Y}^k(Y)       \ar[r]^{(\pi|_E)^*}& H_{\tilde{\theta}|_E}^k(E)   \ar[r]^{} &  \textrm{coker} (\pi|_E)^*    \ar[r]& 0. }
\end{equation}
Notice that $\mathcal{O}_E(-1)=\mathcal{O}_{\widetilde{X}}(E)|_E$, so $h=[E]|_E$.
By Lemma \ref{key}, $i_E^*i_{E*}(\bullet)=h\cup\bullet$ on $H_{\tilde{\theta}|_E}^*(E)$.
Suppose $\pi^*\alpha_k+\sum\limits_{i=1}^{r-1}i_{E*}\left(h^{i-1}\cup(\pi|_E)^*\beta_{k-2i}\right)=0$, where $\alpha_k\in H_{\theta}^k(X)$ and $\beta_{k-2i}\in H_{\theta|_Y}^{k-2i}(Y)$ for $1\leq i\leq r-1$.
Pull it back by $i_E^*$, we get
\begin{displaymath}
(\pi|_E)^*i_Y^*\alpha_k+\sum_{i=1}^{r-1}h^{i}\cup(\pi|_E)^*\beta_{k-2i}=0.
\end{displaymath}
By Corollary \ref{3}, $\beta_{k-2i}=0$ for every $i$.
So $\pi^*\alpha_k=0$. By Corollary \ref{i-s}, $\alpha_k=0$.
Hence  $($\ref{b-u-m}$)$ is injective.
For any $\gamma\in H_{\tilde{\theta}}^k(\widetilde{X})$, by Corollary \ref{3}, there exist $\beta_{k-2i}\in H_{\theta|_Y}^{k-2i}(Y)$ for $0\leq i\leq r-1$, such that $i_E^*\gamma=\sum\limits_{i=0}^{r-1}h^i\cup(\pi|_E)^*\beta_{k-2i}$.
Then
\begin{displaymath}
i_E^*\left[\gamma-\sum_{i=1}^{r-1}i_{E*}\left(h^{i-1}\cup(\pi|_E)^*\beta_{k-2i}\right)\right]=(\pi|_E)^*\beta_k,
\end{displaymath}
which is zero in $\textrm{coker}(\pi|_E)^*$.
From (\ref{commutative1}),
\begin{displaymath}
\gamma-\sum_{i=1}^{r-1}i_{E*}\left(h^{i-1}\cup(\pi|_E)^*\beta_{k-2i})\right)=\pi^*\alpha_k
\end{displaymath}
for some $\alpha_k\in H_{\theta}^k(X)$, which implies  that (\ref{b-u-m}) is surjective.
Hence (\ref{b-u-m}) gives the isomorphism (\ref{isomorphism 1}).

By Proposition \ref{cpb2} and \cite[p. 186, 7.8]{I}, we have the commutative diagram of exact sequences
\begin{displaymath}
\xymatrix{
 \cdots\ar[r]&H_{\theta,c}^k(U)\ar[d]^{\cong} \ar[r]& H_{\theta,c}^k(X) \ar[d]^{\pi^*} \ar[r]& H_{\theta|_Y,c}^k(Y) \ar[d]^{(\pi\mid_E)^*}\ar[r]& H_{\theta,c}^{k+1}(U)\ar[d]^{\cong}\ar[r]&\cdots\\
 \cdots\ar[r]&H_{\tilde{\theta},c}^k(\widetilde{U})       \ar[r]& H_{\tilde{\theta},c}^k(\widetilde{X})     \ar[r]^{i_E^*} & H_{\tilde{\theta}|_E,c}^k(E)     \ar[r]& H_{\tilde{\theta},c}^{k+1}(\widetilde{U})         \ar[r]&\cdots. }
\end{displaymath}
By Corollary \ref{i-s}, $\pi^*$ is injective. By the snake lemma, $i_E^*$ induces an isomorphism
\begin{displaymath}
\textrm{coker}\pi^*\tilde{\rightarrow}\textrm{coker}(\pi|_E)^*.
\end{displaymath}
We get a commutative diagram of short exact sequences
\begin{equation}\label{commutative2}
\xymatrix{
 0\ar[r]&H_{\theta,c}^k(X)\ar[d]^{i_Y^*} \ar[r]^{\pi^*}& H_{\tilde{\theta},c}^k(\widetilde{X})\ar[d]^{i_E^*} \ar[r]& \textrm{coker}\pi^* \ar[d]^{\cong}\ar[r]& 0\\
 0\ar[r]&H_{\theta|_Y,c}^k(Y)       \ar[r]^{(\pi|_E)^*}& H_{\tilde{\theta}|_E,c}^k(E)   \ar[r]^{} &  \textrm{coker} (\pi|_E)^*    \ar[r]& 0. }
\end{equation}
With the same arguments as above, (\ref{b-u-m}) gives the isomorphism (\ref{isomorphism 2}).
\end{proof}

\begin{rem}
After we finished the earlier version \cite{Me1} of the present paper,
Y. Zou \cite{Z} gave the following expression of the blow-up formula on compact complex manifolds with relative cohomology theory
\begin{displaymath}
H_{\tilde{\theta}}^k(\widetilde{X})\tilde{\rightarrow} H_{\theta}^k(X)\oplus \bigoplus_{i=1}^{r-1}H_{\theta|_Y}^{k-2i}(Y)
\end{displaymath}
\begin{equation}\label{rep0}
\quad\quad\quad\alpha\mapsto(\pi_*\alpha,\alpha^{k-2},\ldots,\alpha^{k-2r+2}).
\end{equation}
where   $\alpha^{k-2i}\in H_{\theta|_Y}^{k-2i}(Y)$ for  $0\leq i\leq r-1$ satisfy
\begin{displaymath}
i_E^*\alpha =\sum_{i=0}^{r-1}h^i\cup (\pi|_E)^*\alpha^{k-2i}
\end{displaymath}
for $\alpha\in H_{\tilde{\theta}}^{k}(\widetilde{X})$.
The compactness is necessary there, since the finiteness of the cohomology groups is used in his proof.
Actually, (\ref{rep0}) is inverse to ours  from the proof of  Theorem \ref{1.3}, so it is still an isomorphism on an arbitrary complex manifold.
\end{rem}

\section{Comparisons of three methods}
Suppose that $X$ is a smooth manifold and  $\mathcal{V}$ is a local system of $\mathbb{R}$-modules of finite rank on $X$.
Then $\mathcal{V}$ has a natural soft resolution
\begin{displaymath}
\begin{aligned}
&\quad\qquad\qquad0\rightarrow\mathcal{V}\rightarrow (\mathcal{V}\otimes \mathcal{A}_X^\bullet,id\otimes d)\\
&\mbox{ (or also }0\rightarrow\mathcal{V}\rightarrow (\mathcal{V}\otimes \mathcal{D}_X^{\prime\bullet},id\otimes d), \mbox{  if } X \mbox{  is oriented)}.
\end{aligned}
\end{displaymath}
Set $T^*(X,\mathcal{V})=H^*((\mathcal{V}\otimes \mathcal{A}_X^\bullet,id\otimes d))$.
Then $H^*(X,\mathcal{V})\cong T^*(X,\mathcal{V})$.

Let $\theta$ be a closed one-form on $X$.
The weight $\theta$-sheaf $\underline{\mathbb{R}}_{X,\theta}$ is a local system of $\mathbb{R}$-modules of  rank $1$ on $X$.
We can compute $H^*(X,\underline{\mathbb{R}}_{X,\theta})$ by $T^*(X,\underline{\mathbb{R}}_{X,\theta})$ and also by $H^*_{\theta}(X)$.
Clearly, $f\otimes\alpha\mapsto f\cdot\alpha$ for all $f\in \underline{\mathbb{R}}_{X,\theta}$, $\alpha\in\mathcal{A}^{\bullet}_X$
(or also $\alpha\in\mathcal{D}^{\prime\bullet}_X$,  if $X$ is oriented) give a morphism
\begin{displaymath}
\begin{aligned}
&\quad\qquad\qquad (\underline{\mathbb{R}}_{X,\theta}\otimes \mathcal{A}_X^\bullet,id\otimes d)\rightarrow (\mathcal{A}_X^\bullet, d_\theta)\\
&\mbox{ (or also }(\underline{\mathbb{R}}_{X,\theta}\otimes \mathcal{D}_X^{\prime\bullet},id\otimes d)\rightarrow (\mathcal{D}_X^{\prime\bullet}, d_\theta), \mbox{  if } X \mbox{  is oriented)}
\end{aligned}
\end{displaymath}
between the two soft resolutions of $\underline{\mathbb{R}}_{X,\theta}$, which induces an isomorphism $T^*(X,\underline{\mathbb{R}}_{X,\theta})\tilde{\rightarrow} H^*_{\theta}(X)$.

Now, there are three cohomology groups: $H^*_{\theta}(X)$, $T^*(X,\mathcal{V})$ and $H^*(X,\mathcal{V})$,
which are all the same up to isomorphisms for $\mathcal{V}=\underline{\mathbb{R}}_{X,\theta}$.
In the present work and the subsequent works \cite{Me2, CY, Me3},  the following  two combinations of cohomology groups

$(1)$ $H^*_{\theta}(X)$ and $H^*(X,\mathcal{V})$, where  $\mathcal{V}=\underline{\mathbb{R}}_{X,\theta}$, and

$(2)$ $T^*(X,\mathcal{V})$ and $H^*(X,\mathcal{V})$.\\
were studied respectively, where we can obtain the fruitful results on $H^*(X,\mathcal{V})$ by the sheaf theory.
Evidently, the results obtained by the latter apply to more general cases and coincide with the ones obtained by the former  when we take $\mathcal{V}=\underline{\mathbb{R}}_{X,\theta}$.
Despite all this, the study on  $H^*_{\theta}(X)$ itself can't be replaced by the ones on $T^*(X,\mathcal{V})$ and $H^*(X,\mathcal{V})$.
Many properties of $H^*_{\theta}(X)$ is difficult to  obtain by studying $T^*(X,\mathcal{V})$ and $H^*(X,\mathcal{V})$,
because some subtle techniques only existing among smooth forms can be used on $H^*_{\theta}(X)$ but may not apply to  $T^*(X,\mathcal{V})$  and $H^*(X,\mathcal{V})$ for a general local system $\mathcal{V}$.
This point can be seen from the proofs of Theorem \ref{L-H} and Lemmas \ref{key}, \ref{1}.
Taking Theorem \ref{L-H} for example, we look back on its proof:
To check the commutative diagram, we eliminated the effect of the twisted part coming from the total space $E$ (i.e., $\Omega$)  by the relations (\ref{subtle-1})-(\ref{subtle-3}), which mixed the information of twisted parts of both the base space $X$ and the total space $E$ (i.e., $\theta$ and $\Omega$).
Notice that, a local representative of a class of $H^*_{\theta}(X)$ is just a smooth form while the one of a class of $T^*(X,\mathcal{V})$ is a smooth form twisted with a section of $\mathcal{V}$.
There seems not to be the analogous  relations with (\ref{subtle-1})-(\ref{subtle-3}) for smooth forms twisted with  sections of sheaves.
At this point, the sheaf theory seems also not to work.
Therefore, Theorem \ref{L-H} is difficult to  be completely generalized to $T^*(X,\mathcal{V})$ or $H^*(X,\mathcal{V})$.
Of course, if only the case with a trivial twisted part on $E$ was considered, we will get a partial generalization of Theorem \ref{L-H}; see \cite[Theorem 5.6]{Me2}\cite[Theorem 5.2]{Me3}.
As we saw  above, it is still necessary to study  $H^*_{\theta}(X)$ itself.

\subsection*{Acknowledgements}
The author would like to express his gratitude to the  referees for their helpful comments and suggestions which greatly improved the manuscript.
The author is supported by the National Natural Science Foundation of China (Grant No. 12001500, 12071444), the Scientific and Technological Innovation Programs of Higher Education Institutions in Shanxi (Grant No. 2020L0290) and the Fundamental Research Program of Shanxi Province (Grant No. 201901D111141).


\end{document}